\def\I{i}
\def\etc{\hbox{\rm etc.}}
\def\atan{\hbox{\rm Atan}}
\def\xzo{$x=0$}
\def\xon{$x=1$}
\def\xoo{$x=\infty$}
\def\xpmp{x^p+x^{-p}}				% Special numerator
\def\torzp{x^p-2\cos\zeta+x^{-p}}		% General numerator
\def\tortn{x^n-2\cos\theta+x^{-n}}		% Denominator
\def\tornn{x^{2n}-2x^n\cos\theta+1}		% Denominator * x^n
\def\torone{x^1-2\cos\omega+x^{-1}}		% Factor of denominator
\def\cl#1{\hfill{\bf #1}\hfill}		% \centerline for Plain teX
\def\tfr#1/#2{{#1}/{#2}}			% Text fraction
\def\fr#1/#2{{{#1}\over{#2}}}			% Fraction
\def\intofx#1/#2{\int \fr {dx}/x \, \fr {#1}/{#2}}
\def\intofy#1/#2{\int \fr {dy}/y \, \fr {#1}/{#2}}
\def\intofz#1/#2{\int \fr {dz}/z \, \fr {#1}/{#2}}

\def\PN{\par\noindent}
\def\LP{\par\medskip\noindent\S}

\def\fourier#1{\qquad\qquad\left\{\,#1\,\right\}}
\let\sysfootnote=\footnote
\newcount\footnotenumber
\def\footnote#1{\global\advance\footnotenumber by 1\sysfootnote{$^{\the\footnotenumber}$}{\rm#1}}

\def\Latex#1{}
\Latex{\documentclass{article}}
\Latex{\pagestyle{empty}}

\Latex{\def\cl#1{\centerline{\bf #1}}}
\Latex{\let\footnote=\sysfootnote}
\Latex{\let\nopagenumbers=\relax}
\nopagenumbers

\Latex{
\begin{document}}

% 6.		Choose output in Latin, French, English, or any combination, including none,
%		by defining language macros appropriately (the setup here is for English).
%		For Latin, disable the \fourier displays and cut off after Appendix A.
%
\def\L#1/L{}
\def\F#1/F{}
\def\E#1/E{\PN #1}			% English is the default
%English: \def\E#1/E{\PN #1}
%Latin: \def\L#1/L{\PN #1}
%Latin: \def\E#1/E{}
%Francais: \def\F#1/F{\PN #1}
%Francais: \def\E#1/E{}

% 7.		Support for 11 accented characters in Plain TeX or LaTeX for the French version
%
%ISO-8859: \catcode228=\active  \def^^e4{\"a}    %ISO-8859
%ISO-8859: \catcode224=\active  \def^^e0{\`a}    %ISO-8859
%ISO-8859: \catcode231=\active  \def^^e7{\c{c}}  %ISO-8859
%ISO-8859: \catcode232=\active  \def^^e8{\`e}    %ISO-8859
%ISO-8859: \catcode233=\active  \def^^e9{\'e}    %ISO-8859
%ISO-8859: \catcode234=\active  \def^^ea{\^e}    %ISO-8859
%ISO-8859: \catcode244=\active  \def^^f4{\^o}    %ISO-8859
%ISO-8859: \catcode249=\active  \def^^f9{\`u}    %ISO-8859
%ISO-8859: \catcode251=\active  \def^^fb{\^u}    %ISO-8859
%ISO-8859: \catcode238=\active  \def^^ee{\^\i{}} %ISO-8859
%ISO-8859: \catcode176=\active  \def^^b0{^\circ} %ISO-8859

%IBM-850: \catcode130=\active  \def^^82{\'e}    %IBM-850
%IBM-850: \catcode132=\active  \def^^84{\"a}    %IBM-850
%IBM-850: \catcode133=\active  \def^^85{\`a}    %IBM-850
%IBM-850: \catcode135=\active  \def^^87{\c{c}}  %IBM-850
%IBM-850: \catcode136=\active  \def^^88{\^e}    %IBM-850
%IBM-850: \catcode138=\active  \def^^8a{\`e}    %IBM-850
%IBM-850: \catcode147=\active  \def^^93{\^o}    %IBM-850
%IBM-850: \catcode151=\active  \def^^97{\`u}    %IBM-850
%IBM-850: \catcode150=\active  \def^^96{\^u}    %IBM-850
%IBM-850: \catcode140=\active  \def^^8c{\^\i{}} %IBM-850
%IBM-850: \catcode248=\active  \def^^f8{^\circ} %IBM-850

%=======================================================================================
%
\L{\cl{	Methodus facilis inveniendi integrale huius formulae		}}/L
\F{\cl{	M\'ethode facile pour trouver l'int\'egrale de la formule	}}/F
\E{\cl{	An easy method for finding the integral of the formula	}}/E
$$
    \intofx { x^{n+p}-2x^n\cos\zeta+x^{n-p} }/{ x^{2n}-2x^n\cos\theta+1 },
 $$
\L{\cl{	casu quo post integrationem ponitur vel {\xon} vel \xoo.	}}/L
\F{\cl{	avec borne sup\'erieure d'int\'egration {\xon} ou \xoo.	}}/F
\E{\cl{	with upper limit of integration {\xon} or \xoo.			}}/E
\par\medskip\noindent
\L{\cl{	Auctore L. Eulero.	}}/L
\F{\cl{	Auteur L. Euler.		}}/F
\E{\cl{	Author L. Euler.		}}/E
\par\medskip\noindent
\L{\cl{	Nova acta academiae scientarum Petropolitanae 3 (1785), 1788, p. 3--24.}}/L
\F{\cl{	Nova acta academiae scientarum Petropolitanae 3 (1785), 1788, p. 3--24.}}/F
\E{\cl{	Nova acta academiae scientarum Petropolitanae 3 (1785), 1788, p. 3--24.}}/E
\par\medskip\noindent
\L{\cl{	Conventui exhibita die 18 Martii 1776.	}}/L
\F{\cl{	Pr\'esent\'ee le 18 Mars 1776.		}}/F
\E{\cl{	Read on 18 March 1776.				}}/E
\par\medskip\noindent
\L{\cl{	Scriptor J. G\'elinas, Februarius 2012.		}}/L
\F{\cl{	Traduction par J. G\'elinas, f\'evrier 2012.	}}/F
\E{\cl{	Translation by J. G\'elinas, February 2012.	}}/E
\par\medskip\noindent
%
%====================================================  (Page 3)
%
%---------------------------------------  (Section 1)
%
\L{\LP 1. Denotet $S$ integrale huius formulae generaliter sumtum, ita ut quaeri debeat valor ipsius $S$, casu quo statuitur {\xon}; ubi primum observo, formam propositam multo concinniorem reddi si fractionis numerator et denominator per $x^n$ dividantur; tum enim habebimus}/L
\F{\LP 1. D\'esignant par $S$ l'int\'egrale ind\'efinie de cette formule, l'on peut se demander quelle devrait \^etre la valeur de $S$, dans le cas o\`u {\xon}; notons tout d'abord que la formule propos\'ee se simplifie grandement si le num\'erateur et le d\'enominateur de la fraction sont divis\'es par $x^n$; nous aurons alors\footnote{Voir l'annexe D pour les formules entre accolades~: $a=\pi-\theta, c=\pi-\zeta$ ici -- traducteur.}}/F
\E{\LP 1. Denoting by $S$ the indefinite integral of this formula, we can ask what will be the value of $S$, in the case when {\xon}; first observe that the given formula is much simplified if the numerator and the denominator of the fraction are divided by $x^n$; we will then have\footnote{Refer to Appendix D for the formulas inside curly brackets~: $a=\pi-\theta, c=\pi-\zeta$ here -- translator.}}/E
$$
	S = \intofx { \torzp }/{ \tortn }.
	\fourier{\int_0^\infty \fr { \cosh pt + \cos c }/{ \cosh nt + \cos a } \, dt}
 $$
\L{Hic statim patet, denominatorem $\tortn$ semper in $n$ factores resolui posse, qui singuli sint formae $\torone$, ubi angulum $\omega$ ita capi oportet, ut, dum iste evanescit, simul etiam ipse denominator ad nihilum redigatur.}/L
\F{Il est imm\'ediatement \'evident que le d\'enominateur, $\tortn$, peut toujours \^etre d\'ecompos\'e en $n$ facteurs, qui ont chacun la forme $\torone$, o\`u l'angle $\omega$ doit \^etre choisi de telle mani\`ere que, si cette expression s'annule, le d\'enominateur est r\'eduit en m\^eme temps \`a z\'ero.}/F
\E{It is immediately obvious that the denominator, $\tortn$, can always be factored into $n$ factors, each having the form $\torone$, where the angle $\omega$ must be chosen so that, if this expression is zero, the denominator is at the same time reduced to zero.}/E
%
%====================================================  (Page 4)
%
%---------------------------------------  (Section 2)
%
\L{\LP 2. Posito autem isto factore $\torone=0$, inde sit $x=\cos{\omega}+\I\sin{\omega}$, unde in genere colligitur $x^{ \lambda}=\cos{\lambda\omega}+\I\sin{\lambda\omega}$ et $x^{-\lambda}=\cos{\lambda\omega}-\I\sin{\lambda\omega}$. Hinc ergo denominator istum accipiet valorem: $2\cos n\omega - 2\cos\theta$, qui igitur evanescet, si pro $n\omega$ sumatur aliquis es his valoribus~:}/L
\F{\LP 2. Si toutefois le facteur $\torone=0$, et donc si $x=\cos{\omega}+\I\sin{\omega}$, nous obtenons en g\'en\'eral $x^{ \lambda}=\cos{\lambda\omega}+\I\sin{\lambda\omega}$ et $x^{-\lambda}=\cos{\lambda\omega}-\I\sin{\lambda\omega}$. Ainsi le d\'enominateur aura cette valeur~: $2\cos n\omega - 2\cos\theta$, qui par cons\'equent s'annule, si nous prenons pour $n\omega$ l'une des valeurs~:}/F
\E{\LP 2. In any case if the factor $\torone=0$, and thus if $x=\cos{\omega}+\I\sin{\omega}$, we obtain in general $x^{\lambda}=\cos{\lambda\omega}+\I\sin{\lambda\omega}$ and $x^{-\lambda}=\cos{\lambda\omega}-\I\sin{\lambda\omega}$. Thus the denominator will have this value~: $2\cos n\omega - 2\cos\theta$, which consequently becomes zero, when we choose for $n\omega$ one of the values~:}/E
$$
	\theta, \theta+2\pi, \theta+4\pi, \theta+6\pi, \theta+8\pi, \etc
 $$
\L{quare cum numerus horum valorum debeat esse $=n$, omnes valores anguli $\omega$ erunt sequentes~:}/L
\F{donc comme le nombre de ces valeurs doit \^etre $n$, les valeurs possibles de l'angle $\omega$ seront les suivantes~:}/F
\E{and as the number of these values must be $n$, all the possible values of the angle $\omega$ will be the following~:}/E
$$
	\fr \theta/n, \fr { \theta+2\pi }/n, \fr { \theta+4\pi }/n, \ldots, \fr { \theta+2(n-1)\pi }/n.
 $$
\L{Praeterea vero cum sit $\cos n\omega=\cos\theta$, erit etiam $\sin n\omega=\sin\theta$.}/L
\F{En outre quand $\cos n\omega=\cos\theta$, nous aurons aussi $\sin n\omega=\sin\theta$.}/F
\E{Moreover when $\cos n\omega=\cos\theta$, we will also have $\sin n\omega=\sin\theta$.}/E
%
%---------------------------------------  (Section 3)
%
\L{\LP 3. Cum igitur denominator habeat $n$ factores huius formae~: $\torone$, nostra fractio, quicunque fuerit eius numerator, in $n$ fractiones simplices resolvi poterit, quarum denominatores erunt illi $n$ factores denominatoris. Scribamus igitur brevitatis gratia $\Pi$ loco numeratoris $\torzp$, atque haec fractio~:}/L
\F{\LP 3. Par cons\'equent, puisque le d\'enominateur a $n$ facteurs de la forme $\torone$, notre fraction, quel que soit le num\'erateur, peut se d\'ecomposer en $n$ fractions simples, dont les d\'enominateurs sont les $n$ facteurs du d\'enominateur.  Par cons\'equent en \'ecrivant pour abr\'eger $\Pi$ \`a la place du num\'erateur $\torzp$, cette fraction:}/F
\E{\LP 3. Thus, since the denominator has $n$ factors of the form  $\torone$, our fraction, whatever its numerator, can be expanded into $n$ simple fractions, whose denominators are the $n$ factors of the denominator.  Thus writing for brevity $\Pi$ to replace the numerator $\torzp$, this fraction:}/E
$$
	\fr \Pi/{ \tortn }
 $$
\L{resolvetur in $n$ fractiones simplices, quarum singulae hanc habebunt formam~:}/L
\F{sera d\'ecompos\'ee en $n$ fractions simples, dont chacune aura cette forme~:}/F
\E{will be decomposed into $n$ simple fractions, each of which having this form~:}/E
$$
	\fr P/{ \torone },
 $$
\L{quocirca statuamus~:}/L
\F{et on a donc \'etabli~:}/F
\E{and therefore we have established~:}/E
$$
	\fr \Pi/{ \tortn } = \fr P/{ \torone } + R ;
 $$
\L{ubi littera $R$ omnes reliquas complectatur fractiones, unde statim habebimus}/L
\F{o\`u la lettre $R$ regroupe toutes les autres fractions, et ainsi nous aurons}/F
\E{where the letter $R$ regroups all the other fractions, and thus we shall have}/E
$$
	\fr { \Pi\,(\torone) }/{ \tortn } = P + R\,(\torone).
 $$
%
%====================================================  (Page 5)
%
\L{Quodsi iam faciamus }/L
\F{Si maintenant nous faisons }/F
\E{If now we make }/E
$$
	\torone = 0,
 $$
\L{hinc colligemus numeratorem $P$~; erit enim}/L
\F{nous trouvons le num\'erateur $P$~; car il sera}/F
\E{we find the numerator $P$~; for it will be}/E
$$
	P = \fr { \Pi\,(\torone) }/{ \tortn },
 $$
\L{siquidem in hac aequatione ponatur}/L
\F{quand nous posons dans cette \'equation}/F
\E{when we make in this equation}/E
$$
	x = \cos\omega + \I\sin\omega.
 $$
%
%---------------------------------------  (Section 4)
%
\L{\LP 4. At vero iam vidimus, si ipsi $x$ hunc valorem tribuamus, illius fractionis tam numeratorem quam denominatorem evanescere, quamobrem secundum regulam notissimam loco numeratoris et denominatoris sua scribamus differentialia, ac prodibit}/L
\F{\LP 4. Mais, comme nous l'avons vu, si on leur attribue cette valeur $x$, tant le num\'erateur que le d\'enominateur de cette fraction s'annuleront, et, selon la r\`egle bien connue\footnote{R\`egle de L'H\^opital dans Inst. Calculi Differentialis II -- E212, 1755, ch. 15. Voir: William Dunham, When Euler Met l'H\^opital, Mathematics Magazine, v. 82, n. 1, f\'evrier 2009, pp. 16-25 -- traducteur.}, il faut \'ecrire les quotients diff\'erentiels au lieu du num\'erateur et du d\'enominateur, ce qui produit}/F
\E{\LP 4. But, then, as we have seen, if we attribute to them this value $x$, the numerator as well as the denominator of this fraction will be zero, and according to the well known rule\footnote{L'H\^opital's rule in Inst. Calculi Differentialis II -- E212, 1755, ch. 15. See: William Dunham, When Euler Met l'H\^opital, Mathematics Magazine, Volume 82, Number 1, February 2009, pp. 16-25 -- translator.}, we must write the differential quotients instead of the numerator and the denominator, which gives}/E
$$
	P = \fr { \Pi\,( x^1 - x^{-1} ) }/{ nx^{n} - nx^{-n} }.
 $$
\L{Nunc igitur si loco $x$ valor assignatus scribatur, primo pro $\Pi$ nanciscemur hunc valorem~:}/L
\F{Maintenant, donc, si on attribue cette valeur assign\'ee \`a $x$ on obtient d'abord pour $\Pi$ cette valeur~:}/F
\E{And now, if this chosen value is attributed to $x$ we obtain first for $\Pi$ this value~:}/E
$$
	\Pi = 2 \cos p\omega - 2 \cos \zeta ;
 $$
\L{ex fractione autem oritur iste valor~: }/L
\F{et pour l'autre valeur apparaissant dans la fraction pr\'ec\'edente~: }/F
\E{and for the other value appearing in the the preceding fraction~: }/E
$$
	\fr { \sin\omega }/{ n\sin n\omega },
 $$
\L{quae ergo expressio cum sit realis, numerator quaesitus erit}/L
\F{donc, puisque c'est une expression r\'eelle, le num\'erateur cherch\'e sera}/F
\E{thus, since this is a real expression, the required numerator will be}/E
$$
	P = \fr { 2\sin\omega\,(\cos p\omega - \cos\zeta) }/{ n\sin n\omega }.
 $$
\L{Iam autem vidimus esse}/L
\F{Nous avons d\'ej\`a vu que}/F
\E{We have already seen that}/E
$$
	\sin n\omega = \sin\theta
 $$
\L{unde iste numerator erit}/L
\F{et ce num\'erateur sera}/F
\E{and this numerator will be}/E
$$
	P = \fr { 2\sin\omega\,(\cos p\omega - \cos\zeta) }/{ n\sin\theta }.
 $$
%
%---------------------------------------  (Section 5)
%
\L{\LP 5. Quaelibet igitur fractio partialis ex resolutione fractionis propositae oriunda erit huiusmodi~:}/L
\F{\LP 5. C'est pourquoi chaque fraction partielle d\'ecoulant de la d\'ecomposition de la fraction propos\'ee aura cette forme~:}/F
\E{\LP 5. So each partial fraction coming from the expansion of the proposed fraction will have this form~:}/E
$$
	\fr { 2(\cos p\omega - \cos\zeta) }/{ n\sin\theta } \, \fr { \sin\omega }/{ \torone }
 $$
\L{in qua forma si angulo $\omega$ successive tribuantur omnes valores supra assignati, qui errant}/L
\F{de laquelle forme si l'angle $\omega$ prend toutes les valeurs assign\'ees plus haut, qui sont}/F
\E{from which form if the angle $\omega$ is given all the successive values found above, which are}/E
$$
	\fr \theta/n, \fr { \theta+2\pi }/n, \fr { \theta+4\pi }/n, \ldots, \fr { \theta+2(n-1)\pi }/n,
 $$
\L{orientur omnes fractiones partiales, quae in unam summam}/L
\F{appara\^\i{}tront toutes les fractions partielles, qui, rassembl\'ees dans une somme,}/F
\E{will appear all the partial fractions, which collected in one sum}/E
%
%====================================================  (Page 6)
%
\L{collectae ipsam formam propositam}/L
\F{devront produire la forme propos\'ee}/F
\E{must produce the proposed form}/E
$$
	\fr { \torzp }/{ \tortn }
 $$
\L{producere debebunt, unde etiam singulae in $dx/x$ ductae et integratae, tum vero in unam summam collectae, exibebunt integrale quaesitum $S$.}/L
\F{et, qui, multipli\'ees chacune par $dx/x$ et int\'egr\'ees, puis r\'eunies dans une seule somme, fourniront l'int\'egrale cherch\'ee $S$.}/F
\E{and, which each multiplied by $dx/x$ and integrated, then collected in a single sum, will produce the requested integral $S$.}/E
%
%---------------------------------------  (Section 6)
%
\L{\LP 6. Consideremus igitur fractionem~:}/L
\F{\LP 6. Par cons\'equent consid\'erons cette fraction~:}/F
\E{\LP 6. So let us consider this fraction~:}/E
$$
	\fr { \sin\omega }/{ \torone }
 $$
\L{quae ducta in $dx/x$ praebet}/L
\F{qui, multipli\'ee par $dx/x$, donne}/F
\E{which multiplied by $dx/x$ gives}/E
$$
	\fr { \sin\omega\,dx }/{ x^2 - 2x\cos\omega + 1 },
 $$
\L{cuius integrale, ita sumtum ut evanescat posito {\xzo}, constat esse}/L
\F{dont l'int\'egrale, prise pour s'annuler en {\xzo}, est certainement}/F
\E{whose integral, taken so that it vanishes at {\xzo}, is certainly}/E
$$
	\atan \left( \fr { x\sin\omega }/{ 1 - x\cos\omega } \right).
 $$
\L{Hinc igitur ex hac fractione partiali oritur ista pars integralis~:}/L
\F{Ainsi donc de cette fraction partielle vient cette partie de l'int\'egrale~:}/F
\E{Thus from this partial fraction comes this part of the integral~:}/E
$$
	\fr { 2(\cos p\omega - \cos\zeta) }/{ n\sin\theta }
		  \, \atan\left( \fr { x\sin\omega }/{ 1 - x\cos\omega } \right),
 $$
\L{unde ergo facile deducuntur omnes $n$ partes integralis quaesiti, si loco $\omega$ ordine omnes eius valores assignati substituantur atque in unam summam colligantur.}/L
\F{et ensuite toutes les $n$ parties de l'int\'egrale \'etudi\'ee peuvent \^etre d\'eduites, si sont substitu\'ees \`a la place de $\omega$ toutes les valeurs assign\'ees avant d'\^etre rassembl\'ees dans une m\^eme somme.}/F
\E{and now all the $n$ parts of the integral in question can be deduced if all the assigned values are substituted for $\omega$ to be collected in the same sum.}/E
%
%---------------------------------------  (Section 7)
%
\L{\LP 7. Quoniam autem hoc loco eum tantum integralis valorem postulamus, qui oritur posito \xon, hoc casu fiet}/L
\F{\LP 7. Mais parce que l'on cherche ici la valeur de l'int\'egrale, une fois pos\'e {\xon}, ceci produira dans ce cas}/F
\E{\LP 7. But since we are searching for the value of the integral, once substituted {\xon}, this will produce in this case}/E
$$
	\atan\left( \fr{ x\sin\omega }/{ 1-x\cos\omega }\right) = \atan\left( \fr{ \sin\omega }/{ 1-\cos\omega }\right).
 $$
\L{At vero ista formula $\sin\tfr\omega/{(1-\cos\omega)}$ exprimit cotangentem anguli $\tfr\omega/2$, ideoque tangentem anguli $\tfr{(\pi-\omega)}/2$, ita ut casu pars integralis futura sit\footnote{Quae formula non valet, nisi sit $0<\omega<2\pi$, unde sequitur haec conditio $0<\theta<2\pi$. (Alexander Liapounoff, 1920)}}/L
\F{Mais cette formule $\sin\tfr\omega/{(1-\cos\omega)}$ exprime la cotangente de l'angle $\tfr\omega/2$, et donc la tangente de l'angle $\tfr{(\pi-\omega)}/2$, de sorte que cette partie de la future int\'egrale est\footnote{Cette formule n'est pas valide, sauf si $0<\omega<2\pi$, d'o\`u vient cette condition $0<\theta<2\pi$. (Alexander Liapounoff, 1920)}}/F
\E{But this formula $\sin\tfr\omega/{(1-\cos\omega)}$ expresses the cotangent of the angle $\tfr\omega/2$, and thus the tangent of the angle $\tfr{(\pi-\omega)}/2$, so that this part of the future integral is\footnote{This formula is not valid, unless $0<\omega<2\pi$, and thus follows this condition $0<\theta<2\pi$. (Alexander Liapounoff, 1920)}}/E
$$
	\fr { \cos p\,\omega - \cos\zeta }/{ n\sin\theta }\,(\pi - \omega).
 $$
\L{Hic autem in transitu notasse iuvabit, si desideretur integrale pro casu {\xoo}, tum proditurum esse $\atan\left(-\tfr{\sin\omega}/{\cos\omega}\right)$; quia igitur est $-\tfr{\sin\omega}/{\cos\omega}$ tangens anguli $\pi-\omega$, cum ante habuissenus $\tfr{(\pi-\omega)}/2$, hinc patet casu {\xoo} etiam totum integrale duplo maius fore quam casu \xon.}/L
\F{Ici il faut noter en passant que si l'on veut l'int\'egrale pour le cas {\xoo} il faut prendre $\atan\left(-\tfr{\sin\omega}/{\cos\omega}\right)$; par cons\'equent parce que $-\tfr{\sin\omega}/{\cos\omega}$ est la tangente de l'angle $\pi-\omega$, alors qu'on avait avant $\tfr{(\pi-\omega)}/2$, il est donc \'evident que dans le cas {\xoo} l'int\'egrale totale est deux fois plus grande que  celle du cas \xon.}/F
\E{But here it must be noted in passing that, if we want the integral for the case {\xoo}, we must take $\atan\left(-\tfr{\sin\omega}/{\cos\omega}\right)$; thus since $-\tfr{\sin\omega}/{\cos\omega}$ is the tangent of the angle $\pi-\omega$, while we had before $\tfr{(\pi-\omega)}/2$, it is thus evident that in the case {\xoo} the total integral is twice as large as in the case \xon.}/E
%
%====================================================  (Page 7)
%
%---------------------------------------  (Section 8)
%
\L{\LP 8. Tribuamus igitur angulo $\omega$ successive omnes eius valorem, eosque ordine hic sistamus, eritque}/L
\F{\LP 8. Assignons par cons\'equent \`a $\omega$ toutes ses valeurs successives, pr\'esent\'ees ici dans l'ordre, et nous aurons}/F
\E{\LP 8. Let us assign to $\omega$ all its successive values, presented in order here, and we will have}/E
$$
	S = \fr {\cos\tfr{p\theta}/n - \cos\zeta}/{n\sin\theta} \left(\fr{n\pi-\theta}/n \right)
		 + \fr {\cos\tfr{p(\theta+2\pi)}/n - \cos\zeta }/{ n\sin\theta } \left(\fr {(n-2)\pi-\theta}/n \right)
 $$
$$
		 + \fr {\cos\tfr{p(\theta+4\pi)}/n - \cos\zeta}/{n\sin\theta} \left(\fr{(n-4)\pi-\theta}/n \right)
 $$
$$
		 + \fr {\cos\tfr{p(\theta+6\pi)}/n - \cos\zeta}/{n\sin\theta} \left(\fr{(n-6)\pi-\theta}/n \right)
		 + \etc
 $$
\L{quarum formularum numerus debet esse $=n$. Haec autem expressio statim in duas partes distinguatur hoc nodo indicandas~:}/L
\F{o\`u le nombre de termes doit \^etre $n$. Mais cette expression peut se s\'eparer en deux parties ainsi d\'efinies~:}/F
\E{where the number of terms must be $n$. But this expression can be split in two parts, defined as follows~:}/E
$$
	S = \fr Q/{ n\sin\theta } - \fr { R\cos\zeta }/{ n\sin\theta },
 $$
\L{ita ut sit}/L
\F{de sorte que}/F
\E{so that}/E
$$
	Q =     \fr {     n\pi-\theta }/n \cos\fr p/n\theta
		+ \fr { (n-2)\pi-\theta }/n \cos\fr p/n(\theta+2\pi)
 $$
$$
		+ \fr { (n-4)\pi-\theta }/n \cos\fr p/n(\theta+4\pi)
		+ \etc
 $$
$$
	R = \fr { n\pi-\theta }/n + \fr { (n-2)\pi-\theta }/n + \fr { (n-4)\pi-\theta }/n
		+ \fr { (n-6)\pi-\theta }/n + \etc
 $$
\L{ita ut iam nobis incumbatur in valores litterarum $Q$ et $R$ inquirere.}/L
\F{et de sorte que nous devons maintenant nous enqu\'erir des valeurs des lettres $Q$ et $R$.}/F
\E{and so that we must now search for the value of the letters $Q$ and $R$.}/E
%
%---------------------------------------  (Section 9)
%
\L{\LP 9. Primo autem statim patet, valores ipsius $R$ esse progressionem arithmeticam decrescentem differentia $2\pi/n$, unde summa $n$ terminorum erit $=\pi-\theta$, ita ut sit $R=\pi-\theta$. At vero inventio progressionis $Q$ majorem requirit apparatum, quem in finem sequentes investigationes generaliores praemittamus.}/L
\F{\LP 9. Il est imm\'ediatement clair que les valeurs de ce $R$ sont en progression arithm\'etique d\'ecroissante de raison $2\pi/n$, et que la somme de $n$ termes sera $\pi-\theta$, de sorte que $R=\pi-\theta$. Mais la d\'ecouverte de $Q$,  cependant, exige un plus grand d\'eveloppement de formules que nous pr\'esentons pr\'ec\'ed\'ees de recherches g\'en\'erales.}/F
\E{\LP 9. It is immediately clear that the values of this $R$ are in decreasing arithmetic progression with difference $2\pi/n$, and that the sum of $n$ terms will be $\pi-\theta$, so that $R=\pi-\theta$. But the discovery of $Q$, however, needs a greater development of formulas that we present preceded by general research.}/E
\smallskip
%
%====================================================  (Page 8)
%
%---------------------------------------  (Section 10)
%
\L{\LP 10. Consideretur primo progressio ista cosinuum, quorum anguli in progressione arithmetica progrediantur et quorum numerus sit $n$~:}/L
\F{\LP 10. Consid\'erons en premier la progression des cosinus, dont les angles augmentent selon une progression arithm\'etique et dont le nombre est $n$~:}/F
\E{\LP 10. Consider first the progression of cosines, whose angles are increasing in an arithmetic progression and whose number is $n$~:}/E
$$
	t = \cos(\alpha+2\beta) + \cos(\alpha+4\beta) + \cos(\alpha+6\beta) + \ldots + \cos(\alpha+2n\beta).
 $$
\L{Iam multiplicemus utrinque per $2\sin\beta$, et cum sit}/L
\F{Maintenant multiplions des deux cot\'es par $2\sin\beta$, et, comme}/F
\E{Now let us multiply both sides by $2\sin\beta$, and since}/E
$$
	2\sin\beta\cos\gamma = \sin(\beta+\gamma) - \sin(\gamma-\beta),
 $$
\L{proveniet sequens forma~:}/L
\F{on obtient la forme suivante~:}/F
\E{we obtain the following form~:}/E
$$
	2 t \sin\beta = - \sin(\alpha+\beta) + \sin(\alpha+3\beta) + \sin(\alpha+5\beta)
				+ \ldots + \sin(\alpha+(2n+1)\beta)
 $$
$$
				- \sin(\alpha+3\beta) - \sin(\alpha+5\beta) - \ldots
 $$
\L{ubi omnes termini intermedii manifesto se destruunt, ita ut soli extremi remaneant, hinque ergo fiet}/L
\F{o\`u tous les termes interm\'ediaires manifestement se d\'etruisent deux \`a deux, de sorte que seuls les termes extr\^emes demeurent, ainsi on aura donc}/F
\E{where all the intermediate terms cancel, so that only the extreme terms remain, thus we will have}/E
$$
	t = \fr { \sin(\alpha+(2n+1)\beta) - \sin(\alpha+\beta) }/{ 2\sin\beta }.
 $$
%
%---------------------------------------  (Section 11)
%
\L{\LP 11. Deinde vero iidem cosinus combinentur cum numeris naturalibus $1,2,3,4,\ldots,n$, ac statuatur}/L
\F{\LP 11. Puis, cependant, ces cosinus doivent \^etres combin\'es avec les nombres naturels $1,2,3,4,\ldots,n$, pour avoir}/F
\E{\LP 11. Then, moreover, these cosines must be combined with the natural numbers $1,2,3,4,\ldots,n$, to get}/E
$$
	u = 1\cos(\alpha+2\beta) + 2\cos(\alpha+4\beta) + 3\cos(\alpha+6\beta)+ \ldots + n\cos(\alpha+2n\beta)
 $$
\L{qua expressione ducta in $2\sin\beta$, adhibita resolutione qua modo sumus usi, consequemur}/L
\F{cette expression \'etant multipli\'ee par $2\sin\beta$, utilisant le mode de r\'esolution pr\'ec\'edent, on obtiendra}/F
\E{this expression being multiplied by $2\sin\beta$, using the preceding resolution method, we will obtain}/E
$$
	2u\sin\beta = - \sin(\alpha+\beta) + \sin(\alpha+3\beta) + 2\sin(\alpha+5\beta) + \ldots
 $$
$$
			+ \sin(\alpha+(2n+1)\beta) - 2\sin(\alpha+3\beta) - 3\sin(\alpha+5\beta) - \ldots
 $$
\L{quae forma reducitur ad istam~:}/L
\F{cette forme se r\'eduisant \`a celle-ci~:}/F
\E{this form being reduced to this one~:}/E
$$
	n\sin(\alpha+(2n+1)\beta) - 2u\sin\beta = \sin(\alpha+\beta) + \sin(\alpha+3\beta)
								+ \ldots + \sin(\alpha+(2n-1)\beta)
 $$
\L{quae vocetur $=v$.}/L
\F{que nous d\'esignons par $v$.}/F
\E{that we denote by $v$.}/E
%
%====================================================  (Page 9)
%
%---------------------------------------  (Section 12)
%
\L{\LP 12. Nunc ista series denuo ducatur in $2\sin\beta$, et cum in genere sit}/L
\F{\LP 12. Maintenant on multiplie aussi cette s\'erie par $2\sin\beta$, et comme en g\'en\'eral}/F
\E{\LP 12. Now we multiply also this series by $2\sin\beta$, and since in general}/E
$$
	2\sin\beta\sin\gamma = \cos(\gamma-\beta) - \cos(\gamma+\beta)
 $$
\L{nanciscemur}/L
\F{nous trouverons}/F
\E{we will find}/E
$$
	2 v \sin\beta = \cos\alpha - \cos(\alpha+2\beta) - \cos(\alpha+4\beta) - \ldots
				- \cos(\alpha+2n\beta) + \cos(\alpha+2\beta) + \cos(\alpha+4\beta) + \ldots
 $$
\L{unde ob terminos medios omnes se destruentes colligitur}/L
\F{et donc, comme tous les termes interm\'ediaires se d\'etruisent, nous obtenons}/F
\E{and thus since all the intermediate terms cancel each other, we obtain}/E
$$
	v = \fr{ \cos\alpha - \cos(\alpha+2n\beta ) }/{ 2\sin\beta };
 $$
\L{quare cum sit}/L
\F{c'est pourquoi comme}/F
\E{and since}/E
$$
	u = \fr{ n\sin(\alpha+(2n+1)\beta) - v }/{ 2\sin\beta }
 $$
\L{hinc obtinemus}/L
\F{nous avons donc obtenu}/F
\E{we have thus obtained}/E
$$
	u = \fr{ n\sin(\alpha+(2n+1)\beta) }/{ 2\sin\beta }
		- \fr{ \cos\alpha - \cos(\alpha+2n\beta) }/{ 4\sin^2\beta }.
 $$
%
%---------------------------------------  (Section 13)
%
\L{\LP 13. Combinemus nunc ambas summationes modo traditas in genere, ac statuamus}/L
\F{\LP 13. Combinons maintenant les deux sommes de la mani\`ere donn\'ee en g\'en\'eral, et nous \'etablissons}/F
\E{\LP 13. Let us combine the two sums using the general method, and we establish}/E
$$
	V =   (a+ b)\cos(\alpha+2\beta)
	    + (a+2b)\cos(\alpha+4\beta)
	    + (a+3b)\cos(\alpha+6\beta)
	    + \ldots
	    + (a+nb)\cos(\alpha+2n\beta),
 $$
\L{atque evidens est fore $V=at+bu$, unde loco $t$ et $u$ valoribus sustitutis erit}/L
\F{et il est \'evident que $V=at+bu$, et qu'en substituant \`a la place de $t$ et de $u$ leurs valeurs on aura}/F
\E{and it is evident that $V=at+bu$, and that by replacing $t$ and $u$ by their values we will obtain}/E
$$
	V =     \fr {    a\sin(\alpha+(2n+1)\beta)-a\sin(\alpha+\beta) }/{ 2\sin\beta }
		+ \fr { b n \sin(\alpha+(2n+1)\beta) }/{ 2\sin\beta }
		- \fr { b\cos\alpha  -  b\cos(\alpha+2n\beta) }/{ 4\sin^2\beta }
%		- \fr { b\cos\alpha  +  b\cos(\alpha+2n\beta) }/{ 4\sin^2\beta }     in original Latin
 $$
%
%---------------------------------------  (Section 14)
%
\L{\LP 14. Iam satis perspicuum est progressionem, quam supra littera $Q$ exhibuimus, in ista forma generali pro $V$ inventa contineri, quandoquidem utrique idem terminorum numerus $n$ occurrit, atque coefficientes cosinuum seriei $Q$ etiam progressionem arithmeticam constituunt. Quamobrem pro coefficientibus primo faciamus}/L
\F{\LP 14. Maintenant il est bien clair que la progression, que nous avons d\'esign\'ee plus haut par $Q$, est contenue dans cette forme g\'en\'erale $V$, que le m\^eme nombre $n$ de termes apparaissent des deux c\^ot\'es, et que les coefficients de la s\'erie de cosinus $Q$ constituent une progression arithm\'etique. De cette fa\c{c}on faisons d'abord pour les coefficients}/F
\E{\LP 14. Now it is very clear, that the progression that we have denoted by $Q$ above, is contained in the general form $V$, that the same number of terms occurs on both sides, and that the coefficients of the series of cosines $Q$ constitute an arithmetic progression. In this fashion let us make first for the coefficients}/E
$$
	a + b = \fr{ n\pi-\theta }/n
 $$
\L{et}/L
\F{et}/F
\E{and}/E
$$
	a + 2b = \fr{ (n-2)\pi - \theta }/n
 $$
\L{unde deducimus}/L
\F{donc nous d\'eduisons}/F
\E{thus we deduce}/E
%
%====================================================  (Page 10)
%
$$
	b = - \fr { 2\pi }/n
 $$
\L{et}/L
\F{et}/F
\E{and}/E
$$
	a = \fr{ (n+2)\pi-\theta }/n
 $$
\L{Nunc etiam angulos inter se coaequemus faciamusque}/L
\F{Maintenant faisons aussi les angles proportionnels entre eux}/F
\E{Now let us make the angles be proportional to each other}/E
$$
	\alpha + 2\beta = \fr p/n\theta
 $$
\L{et}/L
\F{et}/F
\E{and}/E
$$
	\alpha + 4\beta = \fr p/n(\theta+2\pi)
 $$
\L{unde colligimus $\beta=\pi p/n$, hincque porro}/L
\F{et nous en d\'eduisons que $\beta=\pi p/n$, et donc de plus}/F
\E{and we deduce that $\beta=\pi p/n$, and moreover}/E
$$
	\alpha = \fr p/n(\theta-2\pi) = - \fr p/n(2\pi-\theta)
 $$
\L{hocque modo fiet $V=Q$. At anguli in expressione ipsius $V$ occurentes erunt~: primo}/L
\F{de cette fa\c{c}on on fait $V=Q$. Mais les angles apparaissant dans cette expression $V$ seront~: premi\`erement}/F
\E{in this way we make $V=Q$.  But the angles occurring in this expression $V$ will be~: first}/E
$$
	\alpha + (2n+1)\beta = - \fr p/n(\pi-\theta)
 $$
\L{ubi cum $p$ sit numero integer, postrema pars $2\pi p$, total circumferentiam exprimens, omissa est, ex quo habebimus}/L
\F{o\`u comme $p$ est nombre entier, la derni\`ere partie $2\pi p$  est omise pour exprimer la ciconf\'erence totale, de quoi nous obtenons}/F
\E{where since $p$ is integer number, the last part $2\pi p$ is omitted to express the total circumference, from which we get}/E
$$
	\sin(\alpha+(2n+1)\beta) = - \sin\fr p/n(\pi-\theta).
 $$
\L{Deinde occurit angulus}/L
\F{Ensuite on aura l'angle}/F
\E{Then we have the angle}/E
$$
	\alpha + \beta = - \fr p/n(\pi - \theta)
 $$
\L{cujus sinus est}/L
\F{dont le sinus est}/F
\E{whose sinus is}/E
$$
	\sin (\alpha + \beta) = - \sin\fr p/n(\pi - \theta).
 $$
\L{Denique erit}/L
\F{Enfin on aura}/F
\E{Finally we have}/E
$$
	\alpha + 2n\beta = - \fr p/n(2\pi - \theta),
 $$
\L{et}/L
\F{et}/F
\E{and}/E
$$
	\cos(\alpha + 2n\beta) = \cos\fr p/n(2\pi - \theta).
 $$
\L{His igitur valoribus substitutis prodibit}/L
\F{La substitution de ces valeurs produit}/F
\E{The substitution of these values produces}/E
$$
	Q = V = \fr{ - \fr{(n+2)\pi-\theta}/n \sin \tfr{p(\pi-\theta)}/n
			 + \fr{(n+2)\pi-\theta}/n \sin \tfr{p(\pi-\theta)}/n }/{ \sin \tfr{p\pi}/n }
 $$
$$
		+ \fr{ 2\pi\sin \tfr{p(\pi-\theta)}/n }/{ 2\sin \tfr{p\pi}/n }
		+ \fr{   \fr{2\pi}/n \cos \tfr{p(2\pi-\theta)}/n
			 - \fr{2\pi}/n \cos \tfr{p(2\pi-\theta)}/n }/{ 4\sin^2 \tfr{p\pi}/n },
 $$
%
%====================================================  (Page 11)
%
\L{quae expresso manifesto reducitur ad hanc~:}/L
\F{cette expression se r\'eduisant clairement \`a ceci~:}/F
\E{this expression clearly being reduced to this~:}/E
$$
	Q = V = \fr{ \pi\sin \tfr{p(\pi-\theta)}/n }/{ \sin \tfr{p\pi}/n }.
 $$
%
%---------------------------------------  (Section 15)
%
\L{\LP 15. Inventis igitur valoribus litterarum $Q$ et $R$, valor integralis quem quaerimus pro casu {\xon} erit}/L
\F{\LP 15. Ayant trouv\'e la valeur des lettres $Q$ et $R$, la valeur de l'int\'egrale cherch\'ee\footnote{Voir l'annexe D pour les formules entre accolades~: $b=\tfr p/n$ ici -- traducteur.} sera pour le cas {\xon}}/F
\E{\LP 15. Having found the value of the letters $Q$ and $R$, the value of the integral considered\footnote{Refer to Appendix D for the formulas inside curly brackets~: $b=\tfr p/n$ here -- translator.} for the case {\xon} will be}/E
$$
	S = \fr { \pi\sin\tfr {p(\pi-\theta)}/n }/{ n\sin\theta\sin\tfr {p\pi}/n } 
		- \fr { (\pi-\theta)\cos\zeta }/{ n\sin\theta }.
	\fourier{\fr 1/n \fr { \pi\sin ab }/{ \sin a\sin \pi b } + \fr { \cos c }/n \fr a/{ \sin a }}
 $$
\L{Sin autem integrale quaeratur a termino {\xzo} usque ad \xoo, eius valor duplo major evadet.}/L
\F{Mais si l'int\'egrale cherch\'ee va de {\xzo} jusqu'\`a \xoo, sa valeur sera deux fois plus grande.}/F
\E{But if the requested integral is taken from {\xzo} to \xoo, its value will be twice as large.}/E
%
%---------------------------------------  (Section 16)
%
\L{\LP 16. His iam in genere expeditis, consideremus casum iam saepius tractarum, quo est $\zeta=90^\circ$ et $\theta=90^\circ$, haecque formula integranda proponitur~:}/L
\F{\LP 16. Ces cas g\'en\'eraux trait\'es, consid\'erons maintenant ceux souvent rencontr\'es quand $\zeta=90^\circ$ et $\theta=90^\circ$, de sorte que la formule int\'egrale propos\'ee est}/F
\E{\LP 16. These general cases having been treated, let us consider now those often met when $\zeta=90^\circ$ and $\theta=90^\circ$, so that the proposed integral formula is}/E
$$
	\intofx { \xpmp }/{ x^{n} + x^{-n} },
	\fourier{\int_0^\infty \fr { \cosh pt }/{ \cosh nt } dt}
 $$
\L{atque pro eius valore casus {\xon} habebimus}/L
\F{ainsi nous aurons pour le cas {\xon} la valeur}/F
\E{and thus we will have for the case {\xon} the value}/E
$$
	S = \fr { \pi\sin\tfr {p\pi}/{2n} }/{ n\sin\tfr {p\pi}/n },
 $$
\L{quae ob}/L
\F{qui, parce que}/F
\E{which since}/E
$$
	\sin \fr { p\pi }/n = 2\sin \fr { p\pi }/{ 2n }\cos \fr { p\pi }/{ 2n },
 $$
\L{abit in formulam illam notissimam}/L
\F{donne une formule remarquable}/F
\E{gives us a noteworthy formula}/E
$$
	\fr \pi/{ 2n\cos\tfr { p\pi }/{ 2n } } = \fr \pi/{ 2n } \sec\fr { p\pi }/{ 2n }.
	\fourier{\fr 1/n \, \fr \pi/2 \sec \fr { \pi b }/2 }.
 $$
\L{Sin autem tantum sumamus $\zeta=90^\circ$ ut formula integranda sit}/L
\F{Mais si on prend seulement $\zeta=90^\circ$ pour que la formule \`a int\'egrer soit}/F
\E{But if we only take $\zeta=90^\circ$ so that the formula to be integrated is}/E
$$
	\intofx { \xpmp }/{ \tortn },
	\fourier{\int_0^\infty \fr { \cosh pt }/{ \cosh nt + \cos a } dt}
 $$
\L{eius valor ab {\xzo} usque ad {\xon} extensus erit}/L
\F{sa valeur prise de {\xzo} jusqu'\`a {\xon} sera}/F
\E{its value taken from {\xzo} up to {\xon} will be}/E
$$
	\fr { \pi\sin\tfr {p(\pi-\theta)}/n }/{ n\sin\theta\sin\tfr {p\pi}/n },
	\fourier{\fr 1/n \fr { \pi \sin ab }/{ \sin a\,\sin \pi b }}
 $$
\L{quae expressio reducitur ad hanc~:}/L
\F{cette expression se r\'eduisant \`a ceci~:}/F
\E{this expression being reduced to this~:}/E
$$
	\fr \pi/{ n\sin\theta } 
		\left( \cos\fr { p\theta }/n - \sin\fr { p\theta }/n \cot\fr { p\pi }/n \right).
 $$
%
%====================================================  (Page 12)
\vfill
\eject
\L{\cl{	Observationes in integrationem traditam.	}}/L
\F{\cl{	Observations sur l'int\'egration usuelle.	}}/F
\E{\cl{	Observations on usual integration.		}}/E
\medskip
%
%---------------------------------------  (Section I)
\L{\LP I. Primum hic observo terminum medium in numeratore exhibitum nullo modo integrationem turbare, quoniam, si solus adesset, integratio nulla laboraret difficultatae ; tum enim formula}/L
\F{\LP I. Premi\`erement, notons que le terme du milieu au num\'erateur ne d\'erange aucunement l'int\'egration car s'il \'etait seul l'int\'egration ne poserait pas de difficult\'e; en effet}/F
\E{\LP I. First, let us note that the middle term in the numerator does not inhibit the integration since if it was alone the integration would not pose any problems;  indeed}/E
$$
	\intofx 1/{ \tortn }
	\fourier{\fr 1/2 \int_0^\infty \fr { dt }/{ \cosh n t + \cos a}}
 $$
\L{reducitur ad hanc formam~:}/L
\F{se r\'eduit \`a cette forme~:}/F
\E{is reduced to this form~:}/E
$$
	\int \fr { x^{n-1} dx }/{ x^{2n}-2x^{n}\cos\theta+1 },
	\fourier{\int_0^\infty \fr { e^{-nt} dt }/{ e^{-2nt} + 2 e^{-nt}\cos a + 1}}
 $$
\L{quae posito $x^{n}=y$ abit in hanc~:}/L
\F{qui apr\`es la substitution $x^{n}=y$ devient celle-ci}/F
\E{which after the substitution $x^{n}=y$ becomes this one}/E
$$
	\fr 1/n \int \fr { dy }/{ y^2-2y\cos\theta+1 },
	\fourier{\fr 1/n \int_0^1 \fr { dy }/{ y^2+2y\cos a+1 }}
 $$
\L{cujus integrale est}/L
\F{dont l'int\'egrale est}/F
\E{whose integral is}/E
$$
	\fr 1/{ n\sin\theta } \atan\left(\fr { y\sin\theta }/{ 1 - y\cos\theta }\right),
	\fourier{\fr 1/{ n\sin a } \left.\atan\left(\fr { y\sin a }/{ 1 + y\cos a }\right)\right|_0^1}
 $$
\L{cujus valor, casu {\xon}, sit}/L
\F{dont la valeur, dans le cas {\xon}, sera}/F
\E{whose value, for the case {\xon}, will be}/E
$$
	\fr 1/{n\sin\theta} \atan\left(\fr { \sin\theta }/{ 1-\cos\theta }\right) =
		   \fr { \pi-\theta }/{ 2n\sin\theta }
	\fourier{\fr 1/{2n}\,\fr a/{ \sin a }}
 $$
\L{qui ductus in $-2\cos\zeta$ praebet illam ipsam partem hinc oriundam in forma supra inventa, quamobrem superstuum foret hunc terminum in calculo retinere, unde hanc formam integralem sumus contemplaturi~:}/L
\F{qui multipli\'ee par $-2\cos\zeta$ donne par cons\'equent un r\'esultat de la forme trouv\'ee ci-dessus, de sorte que le terme obtenu sera retenu dans le calcul et nous sommes amen\'es \`a regarder cette int\'egrale~:}/F
\E{which when multiplied by $-2\cos\zeta$ gives a result of the form found above, thus the resulting term will be kept in the computation, and we are led to look at this integral~:}/E
$$
	\intofx { \xpmp }/{ \tortn }
	\fourier{\int_0^\infty \fr { \cosh pt }/{ \cosh nt + \cos a }\,dt}
 $$
\L{cujus valorem casu {\xon} invenimus}/L
\F{dont la valeur, dans le cas {\xon}, a \'et\'e trouv\'ee \'egale \`a}/F
\E{whose value, for the case {\xon}, we have found to be}/E
$$
	= \fr { \pi\sin\tfr {p(\pi-\theta)}/n }/{ n\sin\theta \sin\tfr { p\pi }/n },
	\fourier{\fr 1/n \, \fr { \pi\sin ab }/{ \sin a\,\sin \pi b}}
 $$
\L{quem brevitas gratia littera $P$ designemus, ita ut sit}/L
\F{que pour abr\'eger nous d\'esignons par la lettre $P$, de sorte que}/F
\E{that for brevity we denote by the letter $P$, so that}/E
$$
	P = \fr { \pi\sin\tfr {p(\pi-\theta)}/n }/{ n\sin\theta\sin\tfr {p\pi}/n }.
 $$
\L{Tum vero etiam invenimus casu {\xoo} valorem huius formulae esse $2P$.}/L
\F{Puis, effectivement, nous avons aussi trouv\'e pour le cas {\xoo} que la valeur de cette formule est $2P$.}/F
\E{And, effectively, we also found that in the case {\xoo} the value of this formula is $2P$.}/E
\medskip
%
%---------------------------------------  (Section II)
%
\L{\LP II. Secundo loco probe notari oportet, exponentem $p$ neccessario minorem esse debere quam exponentem $n$, quia alioquin fractio foret spuria, et variabilis $x$ in numeratore tot vel plures dimensiones esset habitura quam in denominatore. Quoties autem hoc evenit, integrali praeter partes, quas per resolutionem}/L
\F{\LP II. En second lieu, il faut noter que l'exposant $p$ doit n\'ecessairement \^etre plus petit que l'exposant $n$  sinon la fraction serait impropre et la variable $x$ aurait plus de degr\'es dans le num\'erateur que dans le d\'enominateur. Lorsque cela se produit, en plus des termes int\'egr\'es provenant de fractions partielles,}/F
\E{\LP II. Secondly, it must be noted that the exponent $p$ must necessarily be smaller than the exponent $n$ since otherwise the fraction would be improper and the variable $x$ would have more degrees in the numerator than in the denominator. When this happens, in addition to the terms integrated from the partial fractions}/E
%
%====================================================  (Page 13)
%
\par\noindent
\L{in fractiones partiales sumus nacti, quantitas quaedam integra adiici debet, id quod in nostra solutione non est factum, quamobrem tales casus hinc prorsus excludi convenit. Ceterum quolibet casu has partes integras facile erit adiicere ad partes quas nobis nostra methodus suppeditabit.}/L
\F{certaines quantit\'es devraient \^etre ajout\'ees, qui n'\'etaient pas dans notre solution, c'est pourquoi de tels cas doivent \^etre totalement exclus. Mais dans tous les cas, il sera facile d'ajouter ces termes aux parties que notre m\'ethode fournit.}/F
\E{there must be added others, which were not taken into account in our solution, thus this is why such cases must totally excluded. But in all cases, it will be easy to add these terms to the parts provided by our method.}/E
%
%---------------------------------------  (Section III)
%
\L{\LP III. Ex ipsa solutione, quam dedimus, perspicuum est exponentem $p$ necessario integrum statui debere, quia alias operationes ibi exhibitae locum habere non possent; unde eo magis mirum videbitur, quod conclusiones inventae substitere queant, etiamsi iste exponens $p$ fuerit numerus fractus quicunque, dummodo minor quam $n$, propterea quod hos casus semper ad exponentes integros reducere licet. Ad hoc ostendendum ponamus esse $p=\tfr q/\lambda$, atque forma nostra, posito $x=z^{\lambda}$, reducetur ad hanc formam~:}/L
\F{\LP III. Dans cette solution, que nous avons donn\'ee, il est clair que l'exposant $p$ doit n\'ecessairement \^etre un entier parce que sinon les op\'erations effectu\'ees ne pourraient avoir lieu; et il sera vu comme merveilleux que les conclusions trouv\'ees peuvent demeurer, m\^eme quand cet exposant $p$ sera une fraction quelconque, quoique inf\'erieure \`a $n$ parce que ces cas peuvent toujours se r\'eduire au cas d'exposants entiers. Pour le montrer posons $p=\tfr q/\lambda$,  de sorte que notre formule, substituant $x=z^{\lambda}$, est r\'eduite \`a cette forme~:}/F
\E{\LP III. In this solution, that we have given, it is clear that the exponent $p$ must necessarily be an integer since otherwise the operations done could not take place; and it will be seen as marvelous that the conclusions found can stand, even when this exponent $p$ will be any fraction, although smaller that $n$ since these cases can be reduced to the case of integer exponents. To show this let $p=\tfr q/\lambda$, so that our formula, substituting  $x=z^{\lambda}$, will be reduced to this form~:}/E
$$
	\lambda \intofz { z^{q}+z^{-q} }/{ z^{\lambda n}-2\cos\theta+z^{-\lambda n} }
 $$
\L{ubi cum omnes exponentes sint integri, ac pro terminis integrationis, qui erant {\xzo}, {\xon} et \xoo, etiam fiat $z=0, z=1$ et $z=\infty$, pro $z=1$ valor integralis erit}/L
\F{o\`u comme tous les exposants sont entiers, et que les bornes d'int\'egration, qui \'etaient {\xzo}, {\xon} et {\xoo} deviennent alors $z=0, z=1$ et $z=\infty$, la valeur de l'int\'egrale pour $z=1$ sera}/F
\E{where since all the exponents are integers, and the limits of integration, which were {\xzo}, {\xon} and {\xoo} become then $z=0, z=1$ and $z=\infty$, the value of the integral for $z=1$ will be}/E
$$
	\fr{ \lambda\pi\sin\tfr { q(\pi-\theta) }/{ \lambda n } }/{ \lambda n\sin\theta\sin\tfr { q\pi }/{ \lambda n } },
 $$
\L{qui, restituto loco $\tfr q/\lambda$ valore $p$, abit in hunc~:}/L
\F{qui, apr\`es remplacement de $\tfr q/\lambda$ par sa valeur $p$, devient ceci}/F
\E{which, when $\tfr q/\lambda$ is replaced by its value $p$, becomes this}/E
$$
	\fr{ \pi\sin\tfr { p(\pi-\theta) }/n }/{ n\sin\theta\sin\tfr { p\pi }/n },
 $$
\L{quae expressio cum superiore prorsus congruit. Atque hinc intelligitur, quominus etiam exponenti $p$ valores irrationales tribuantur, dumne superent exponentem $n$, semper hoc evenire debere.}/L
\F{une expression qui est tout en accord avec ce qui pr\'ec\`ede. Et donc il est \'egalement entendu que si l'on attribue des valeurs irrationelles \`a l'exposant $p$,  \`a condition qu'elles ne d\'epassent pas l'exposant $n$, cela se produira encore.}/F
\E{an expression that is in total accordance with the preceding. And thus it is also agreed that when irrational values are attributed to the exponent $p$, provided they do not exceed the exponent $n$, this will still happen.}/E
%
%---------------------------------------  (Section IV)
%
\L{\LP IV. Hic iam quaestio oritur maximi momenti, utrum etiam exponenti $p$ dare liceat valores imaginarios nec ne ? Hoc autem affirmandum videtur,  quandoquidem imaginaria certe non sint majora quam $n$; unde concludimus, dummodo valor ipsus $p$ ita capiatur imaginarius, ut ipsa formula differentialis}/L
\F{\LP IV. Ici se pose maintenant une question de grande importance, s'il est ou non licite de donner des valeurs imaginaires \`a l'exposant $p$ ? Mais il sera vu qu'on peut r\'epondre positivement, car les imaginaires ne sont certainement pas plus grands que $n$; et nous concluons, sous condition que la valeur de ce $p$ soit imaginaire, de sorte que cette formule diff\'erentielle}/F
\E{\LP IV. Here must be asked a very important question, is it licit or not to give imaginary values to the exponent $p$ ? But it will be seen that we can answer affirmatively,  since the imaginaries are certainly not greater than $n$; and we conclude, on the condition that the value of this $p$ be taken as imaginary, so that this differential formula}/E
%
%====================================================  (Page 14)
%
\L{maneat realis, tum etiam conclusiones nostras veritati consentaneas esse mansuras\footnote{Quamquam hoc non satis confirmatur est, tamen hic revera ponere licet $p=\I q$, quia integrale valore absoluto ipsius $p$ manente satis parvo in seriam integram secundum ipsius $p$ potestates progredientem resolvi potest. (Alexander Liapounoff, 1920)}. Hoc autem evenit, si statuamus $p=\I q$; tum enim, cum in genere sit}/L
\F{demeure r\'eelle, que nos conclusions restent aussi valides\footnote{Bien que ce ne soit pas suffisamment justifi\'e, il est cependant en fait permis de poser ici $p=\I q$, parce que l'int\'egrale, la valeur absolue de ce $p$ restant assez petite, peut se d\'evelopper en fonction de ce $p$ en s\'erie de puissances enti\`eres successives. (Alexander Liapounoff, 1920)}. Ceci arrive, si on pose $p=\I q$; car alors, comme on a en g\'en\'eral}/F
\E{remains real, that our conclusions also remain valid\footnote{Although this is not sufficiently justified, it is still in fact permitted to put here $p=\I q$, since the integral, if the absolute value of $p$ is kept small enough, can be developed in terms of this $p$ as a series of successive integer powers. (Alexander Liapounoff, 1920)}. This happens, if we set $p=\I q$; since then, as we have in general}/E
$$
	e^{\I\phi} + e^{-\I\phi} = 2\cos\phi,
 $$
\L{quia nostro casu est $\phi=q\log x$, ipsa formula integralis erit}/L
\F{qui dans notre cas est $\phi=q\log x$, la formule int\'egrale elle-m\^eme sera}/F
\E{which in our case is $\phi=q\log x$, the integral formula itself will be}/E
$$
	\intofx { 2\cos (q\log x) }/{ \tortn }.
 $$
\L{Nun igitur videamus quamnam formam nostrum integrale casu {\xon} sit recepturum, et quoniam sinus angulorum imaginarium sunt etiam imaginarii, quandoquidem}/L
\F{Maintenant nous voyons ce que notre formule int\'egrale devient dans le cas {\xon}, et pour les sinus d'angles imaginaires qui sont aussi imaginaires, puisque}/F
\E{Now we see what our integral formula  becomes in the case {\xon}, and for sinus of imaginary angles which are also imaginary, since}/E
$$
	e^{\I\phi} - e^{-\I\phi} = \fr 2/\I \sin\phi,
 $$
\L{loco $\phi$ scribamus $\I\psi$, eritque}/L
\F{\`a la place de $\phi$ \'ecrivons $\I\psi$, de sorte que}/F
\E{instead of $\phi$ let us write $\I\psi$, so that}/E
$$
	\sin\I\psi = \fr \I/2 (e^{-\psi} - e^{+\psi}),
 $$
\L{unde in forma integrali erit}/L
\F{et dans la forme int\'egrale on aura}/F
\E{and in the integral form we will have}/E
$$
	\fr p/n = \I \fr q/n,
 $$
\L{ideoque loco $\psi$ scribamus $\fr q/n(\pi-\theta)$ pro numeratore, et $\fr{q\pi}/n$ pro denominatore, ex quo valor integralis ab {\xzo} ad {\xon} extendus erit}/L
\F{de m\^eme \`a la place  de $\psi$ \'ecrivons $\fr q/n(\pi-\theta)$ au num\'erateur, et $\fr{q\pi}/n$ au d\'enominateur, d'o\`u la valeur de l'int\'egrale prise de {\xzo} jusqu'\`a {\xon} sera}/F
\E{similarly in place of $\psi$ let us write $\fr q/n(\pi-\theta)$ in the numerator, and $\fr{q\pi}/n$ in the denominator, from which the value of the integral taken from {\xzo} up to {\xon} will be}/E
$$
	\fr \pi/{ n\sin\theta }\,
	\fr{ e^{- \fr q/n(\pi-\theta) } - e^{+ \fr q/n(\pi-\theta) } }/
		{ e^{- \fr {q\pi}/n } - e^{+ \fr {q\pi}/n } }.
 $$
\L{Hinc igitur formemus sequens Theorema notaru dignissimum~:}/L
\F{Ainsi donc formulons le th\'eor\`eme remarquable suivant~:}/F
\E{Thus let us formulate the following noteworthy theorem~:}/E
\par\medskip\noindent
\L{\cl{Quodsi ista formula integralis~:}}/L
\F{\cl{Si la formule int\'egrale~:}}/F
\E{\cl{If the integral formula~:}}/E
$$
	\intofx { \cos (q\log x) }/{ \tortn }
	\fourier{ \fr 1/2 \int_0^\infty \fr { \cos qt }/{ \cosh nt + \cos a } \,dt}
 $$
\L{\cl{a termino {\xzo} usque ad {\xon} entendatur, eius valor semper erit}}/L
\F{\cl{est prise de la limite {\xzo} jusqu'\`a {\xon}, sa valeur sera toujours}}/F
\E{\cl{is taken from the limit {\xzo} to {\xon}, its value will always be}}/E
$$
	\fr \pi/{ 2 n \sin\theta }\,
	\fr{ e^{- \fr q/n(\pi-\theta) } - e^{+ \fr q/n(\pi-\theta) } }/
		{ e^{- \fr {q\pi}/n } - e^{+ \fr {q\pi}/n } }.
	\fourier{ \fr 1/{ 2n } \fr{ \pi \sinh(\tfr aq/n)}/{ \sin a \sinh(\pi \tfr q/n)} }
 $$
\L{\cl{Sin autem integrale extendatur ab {\xzo} usque ad \xoo,}}/L
\F{\cl{Mais si l'int\'egrale est prise de {\xzo} jusqu'\`a \xoo,}}/F
\E{\cl{But if the integral is taken from {\xzo} up to \xoo,}}/E
\L{\cl{valor prodibit duplo major.}}/L
\F{\cl{sa valeur sera deux fois plus grande.}}/F
\E{\cl{its value will be twice more.}}/E
\par\medskip\noindent
\L{Hoc theorema utique eo majorem attentionem meretur, quod nulla via patet, eius veritatem directe demonstrandi.}/L
\F{Bien s\^ur, ce th\'eor\`eme m\'erite une plus grande attention, de sorte qu'il ne manque rien \`a la d\'emonstration de sa validit\'e.}/F
\E{Of course, this theorem deserves greater attention, so that nothing would be missing from the demonstration of its validity.}/E
\medskip
%
%====================================================  (Page 15)
%
%---------------------------------------  (Section V)
%
\L{\LP V. Revertamur autem ad formam integralem primo expositam, et quoniam numerator duabus constat partibus $x^p$ et $x^{-p}$, unde summa integralium pro {\xon} inventa est $=P$, at pro casu {\xoo} duplo major $=2P$, hic maxime notatu dignum occurit, quod pro termino {\xoo} utraque pars numeratoris eundem producat valorem $=P$. Semper enim erit, integrale ab {\xzo} ad {\xoo} extendendo,}/L
\F{\LP V. Mais revenons \`a la formule int\'egrale pr\'esent\'ee en premier, dont le num\'erateur a deux termes $x^p$ et $x^{-p}$, et dont la somme int\'egrale pour {\xon} a \'et\'e trouv\'ee $=P$, pour {\xoo} deux fois plus $=2P$, il est particuli\`erement indiqu\'e de noter ici, que pour la borne {\xoo} ces deux parties du num\'erateur produisent une valeur $=P$. On a toujours en effet, prenant l'int\'egrale de {\xzo} \`a {\xoo},}/F
\E{\LP V. But let us revert to the integral form presented first, whose numerator has two terms $x^p$ and $x^{-p}$, and whose integral sum for {\xon} was found to be $=P$, for the case {\xoo} twice more $=2P$, it is particularly noteworthy here, that for the limit {\xoo} these two parts of the numerator produce a value $=P$. We always have in fact, taking the integral from {\xzo} to {\xoo},}/E
$$
	\intofx { x^{+p} }/{ \tortn } = \intofx { x^{-p} }/{ \tortn } = P.
 $$
\L{Ad hoc ostendendum ponamus pro posteriore formula $x=1/z$, eaque induet hanc formam~:}/L
\F{Pour le montrer, posons dans la derni\`ere formule $x=1/z$, et appara\^\i{}tra cette forme~:}/F
\E{To show it, let us put in the lat formula $x=1/z$, and this form appears~:}/E
$$
	- \intofz { z^{+p} }/{ z^{-n}-2\cos\theta+z^{+n} },
 $$
\L{quae cum sit priori formae prorsus similis, solo signo $-$ excepto, eius valor a termino $z=0$ usque ad $z=\infty$, negative sumptus, primae formulae erit aequalis. Cum autem sit $z=1/x$, isti termini integralis erunt ab {\xoo} usque ad {\xzo}, qui ergo si invertantur, etiam signum integralis erit mutandum, sicque ipsi priori formulae aequale evadet; quare cum ambae formulae coniunctae summam habeant $=2P$, utriusque seorsim sumtae valor erit $=P$, unde deducitur sequens theorema notaru pariter dignissimum~:}/L
\F{qui comme elle est tout \`a fait semblable \`a la forme pr\'ec\'edente, except\'e le seul signe $-$, sa valeur de la borne $z=0$ jusqu'\`a $z=\infty$, rendue n\'egative, est \'egale \`a la premi\`ere formule. Mais comme $z=1/x$, ses bornes d'int\'egration \'etant de {\xoo} \`a {\xzo}, qui lorsque intervertis, changent donc le signe de l'int\'egrale, elle est \'egale \`a la formule pr\'ec\'edente; c'est pourquoi comme les deux formules conjugu\'ees ont pour somme $=2P$, chacune a une valeur \'egale \`a $P$, et nous d\'eduisons le th\'eor\`eme suivant m\'eritant aussi d'\^etre connu~:}/F
\E{which as it is totally similar to the previous form, the sign $-$ excepted, has a value from the limit $z=0$ to $z=\infty$, when negated, equal to the first formula. But as $z=1/x$, its limits of integration being from {\xoo} to {\xzo}, which when interchanged, change the sign of the integral, it is equal to the previous formula; that is why as the two formulas conjugated have sum $2P$, each one has a value equal to $P$, and we deduce the following theorem worthy to be also known~:}/E
\medskip
\L{\cl{Istius formulae integralis~:}}/L
\F{\cl{La valeur de cette formule int\'egrale~:}}/F
\E{\cl{The value of this integral formula~:}}/E
$$
	\intofx { x^{\pm p} }/{ \tortn },
	\fourier{\fr 1/2 \, \int_{-\infty}^\infty \fr { e^{\pm pt}\,dt }/{ \cosh nt + \cos a }}
 $$
\L{\cl{valor a termino {\xzo} usque ad {\xoo} entensus semper est}}/L
\F{\cl{prise de la borne {\xzo} jusqu'\`a {\xoo}  est toujours}}/F
\E{\cl{taken from the bound {\xzo} up to {\xoo} is always}}/E
$$
	= P = \fr{ \pi\sin\tfr {p(\pi-\theta)}/n }/{ n\sin\theta\sin \tfr{p\pi}/n }.
	\fourier{\fr 1/n \fr { \pi\sin ab }/{ \sin a\,\sin \pi b }}
 $$
\L{Evidens autem est hanc equalitatem pro casu {\xon}, neutiquam locum habere posse.}/L
\F{Clairement cependant cette \'egalit\'e ne peut pas avoir lieu pour le cas \xon.}/F
\E{Clearly however this equality cannot hold for the case \xon.}/E
\medskip
%
%====================================================  (Page 16)
%
%
%---------------------------------------  (Section VI)
%
\L{\LP VI. Quoniam in nostra formula differentiali tantum occurrit terminus $2\cos\theta$, cuius valor idem manet, etiamsi pro $\theta$ sumeremus $\theta \pm 2k\pi$, maxime hic mirum videri debet, quod tum valor integralis maxime diversus sit proditurus, scilicet}/L
\F{\LP VI. Alors que dans notre formule diff\'erentielle appara\^\i{}t seulement le terme $2\cos\theta$, dont la valeur reste identique, m\^eme si nous prenons $\theta \pm 2k\pi$ \`a la place de $\theta$, il doit sembler bien \'etonnant qu'une valeur bien diff\'erente pour l'int\'egrale soit produite, soit}/F
\E{\LP VI. Whereas in our differential formula there appears only the term $2\cos\theta$, whose value remains identical, when we take $\theta \pm 2k\pi$ to replace $\theta$, it must seem quite surprising that a very different value for the integral is produced, which is}/E
$$
	= \fr{ \pi\sin\tfr {(\pi-\theta\pm 2k\pi)p}/n }/{ n\sin\theta\sin\tfr{\pi p}/n },
 $$
\L{unde merito quaeritur, quisnam horum valorum veritati sit conformis, ad quod certe nihil aliud responderi potest, nisi quod omnes veritati aeque consentanei sint censendi\footnote{Formula inventa nonnisi hac conditione $0<\theta<2\pi$ valet.  (Alexander Liapounoff, 1920)}, id quod eo minus mirum videri debet, quod onmes huiusmodi formulae integrales revera sunt functiones multiformes, atque adeo infinitiformes, id quod ex hoc exemplo simplicissimo~: $\int\fr{dx}/{1+x^2}$ intelligi potest. Cum enim eius integrale exhibeat arcum circuli cuius tangens est $x$, tales autem arcus innumerabiles dentur, quorum eadem sit tangens $=x$, necesse est, ut omnes aeque in hac forma integrali contineantur. Quin etiam in nostro valore invento $P$ loco $\pi$ quoque scribere licet $\pi+2k\pi$, eiusque valor nihilominus cum veritate consistere poterit. Verum in huiusmodi integrationibus perpetuo valores minimi desiderari solent, hocque modo omnis difficultas e medio est sublata.}/L
\F{et il vaut la peine de demander laquelle de ces valeurs est exacte, ce \`a quoi l'on ne peut rien r\'epondre, sauf en v\'erit\'e que toutes doivent \^etre consid\'er\'ees valides\footnote{La formule trouv\'ee est valide seulement si $0<\theta<2\pi$.  (Alexander Liapounoff, 1920)}, et ce qui peut sembler moins \'etrange, que toutes ces formules int\'egrales sont des fonctions multiformes, et donc infini-formes, ce que cet exemple simple~: $\int\fr{dx}/{1+x^2}$ peut aider \`a comprendre. Comme en effet cette int\'egrale produit un arc circulaire dont la tangente est $x$, mais qu'on peut donner d'innombrables arcs dont la tangente est $x$, il est n\'ecessaire que tous soient contenus dans la forme int\'egr\'ee. Mais dans notre valeur trouv\'ee pour $P$ \`a la place de $\pi$ on peut aussi \'ecrire $\pi+2k\pi$, n\'eanmoins sa valeur peut rester valide. Mais dans ce genre d'int\'egrales on d\'esire habituellement les valeurs minimales, de sorte que toute difficult\'e est enlev\'ee de la formule.}/F
\E{and it is worth asking which of these values is exact, to which nothing else can be answered, except the truth that all are to be considered equally valid\footnote{the proposed formula is valid only if $0<\theta<2\pi$.  (Alexander Liapounoff, 1920)}, and which must seem less strange, that all these integral formulas are multiform functions, and thus infiniteform, what this simple example~: $\int\fr{dx}/{1+x^2}$ can help to understand. As in fact this integral produces the circular arc whose tangent is $x$, but as innumerable arcs can be given whose tangent is $x$, it is necessary that all be contained in the integrated form. But in our explicit value for $P$, in place of $\pi$ we can write also write $\pi+2k\pi$, nevertheless its value can remain valid. However in this kind of integrals the minimal values are usually wanted, so that all difficulties are taken out of the formula.}/E
%
%---------------------------------------  (Section VII)
%
\L{\LP VII. Deinde in Analysi supra adhibita supposuimus onmes factores denominatoris inter se esse inaequales, id quod utique semper evenit, nisi sit $\cos\theta=\pm 1$, quippe quibus casibus denominator quadratum involvit~: sit enim is}/L
\F{\LP VII. Ensuite dans l'analyse pr\'esent\'ee plus haut nous avons suppos\'e tous les facteurs du d\'enominateur in\'egaux, ce qui se produit, sauf si $\cos\theta=\pm 1$, dans quels cas le d\'enominateur pr\'esente des carr\'es~: il est alors}/F
\E{\LP VII. In the analysis presented above, we have supposed the factors of the denominator unequal to each other, which happens, unless $\cos\theta=\pm 1$, in which cases the denominator involves squares~: it is then}/E
$$
	= x^{-n}(x^{n}\pm 1)^2;
 $$
\L{ex quo patet omnes factores $x^n\pm 1$ bis occurere debere. Hoc incommodum etiam innuitur per ipsam nostram formulam $P$, }/L
\F{d'o\`u il suit que tous les facteurs $x^n\pm 1$ doivent appara\^\i{}tre deux fois. Ce probl\`eme implique \'egalement pour notre formule $P$}/F
\E{from which follows that all factors $x^n\pm 1$ must appear twice. This problem also implies for our formula $P$}/E
%
%====================================================  (Page 17)
%
\par\noindent
\L{quae casu $\theta=0$ valorem indicat infinitum. Verum posito $\theta=\pi$, singulare phaenomenon se offert, dum formulae pro $P$ inventae tam numerator quam denominator evanescunt, atque adeo fractio determinatum nanciscitur valorem. Ponamus enim $\theta=\pi-\omega$, existente $\omega$ infinite parvo, eritque}/L
\F{que le cas $\theta=0$ indique une valeur infinie. Mais posant $\theta=\pi$, un ph\'enom\`ene singulier se produit, quand dans la formule trouv\'ee pour $P$ tant le num\'erateur que le d\'enominateur s'annulent, de sorte que nous devons d\'eterminer la valeur de leur quotient. Posons en effet $\theta=\pi-\omega$, pour un $\omega$ infiniment petit, et nous aurons}/F
\E{that the case $\theta=0$ signals an infinite value. But putting $\theta=\pi$, a singular phenomenon happens, when in the formula found for $P$ both the numerator and the denominator vanish, so that we must determine the value of their quotient. Let us put to that effect $\theta=\pi-\omega$, for an infinitely small $\omega$, and we will have}/E
$$
	\sin\theta=\sin\omega=\omega;
 $$
\L{at ob $\pi-\theta=\omega$, in numeratore habebimus}/L
\F{et puisque $\pi-\theta=\omega$, on aura au num\'erateur}/F
\E{and since $\pi-\theta=\omega$, we will have in the numerator}/E
$$
	\sin\fr { p\omega }/n = \fr { p\omega }/n,
 $$
\L{unde valor ipsius $P$ resultat}/L
\F{et cette valeur de $P$ r\'esulte}/F
\E{and this value of $P$ results}/E
$$
	\fr {\pi p}/{nn\sin\fr {\pi p}/n},
 $$
\L{qui cum penitus sit determinatus, nullum plane dubium superesse potest, quin cum veritate conspiret, unde sequens enascitur Theorema maxime memorabile~:}/L
\F{qui comme elle est pleinement d\'etermin\'ee, il ne peut rester aucun doute qu'elle soit valide, et que survient le th\'eor\`eme suivant des plus m\'emorables~:}/F
\E{which as it is fully determined, there can remain no doubt that it is valid, and that emerges the following theorem amongst the most memorable~:}/E
\par\medskip\noindent
\L{\cl{ Theorema. Proposita formula differentiali}}/L
\F{\cl{ Th\'eor\`eme. Soit la formule diff\'erentielle}}/F
\E{\cl{ Theorem. Consider the differential formula}}/E
$$
	\fr {dx}/x \, \fr { \xpmp }/{ x^n+2+x^{-n} } = \fr { x^{n-1}\,dx\,(\xpmp) }/{ (x^n+1)^2 },
	\fourier{\int_0^\infty \fr { \cosh(2pt) }/{ \cosh^2(nt) }\,dt}
 $$
\L{\cl{ si eius integrale a termino {\xzo} usque ad {\xon} extendatur, eius valor semper erit}}/L
\F{\cl{ dont l'int\'egrale est prise de {\xzo} jusqu'\`a \xon, sa valeur sera toujours}}/F
\E{\cl{ whose integral is taken from {\xzo} to \xon, its value will always be}}/E
$$
	\fr { \pi \tfr p/n }/{ n \sin \pi \tfr p/n };
	\fourier{\fr 1/n \fr { \pi b }/{ \sin \pi b }}
 $$
\L{\cl{ sin autem usque ad terminum {\xoo} extendatur, eius valor erit duplo major,}}/L
\F{\cl{ mais si prise jusqu'\`a {\xoo}, sa valeur sera deux fois plus grande,}}/F
\E{\cl{ but if taken up to {\xoo}, its value will be twice greater,}}/E
$$
	\fr { 2 \pi \tfr p/n }/{ n \sin \pi \tfr p/n }.
 $$
\medskip\noindent
\L{Demonstratio huius Theorematis directa.}/L
\F{D\'emonstration directe de ce th\'eor\`eme.}/F
\E{Direct demonstration of this theorem.}/E
\par\medskip\noindent
\L{Formula ista integralis resolvatur sequenti modo~:}/L
\F{Cette formule int\'egrale se d\'ecompose ainsi~:}/F
\E{This integral formula can be decomposed thus~:}/E
$$
	\intofx { x^{n+p}+x^{n-p} }/{ (1+x^n)^2 } = \fr Q/{1+x^n} + \intofx R/{1+x^n}.
 $$
\L{Sumantur igitur differentiala simulque ducantur in $\tfr x/{dx}$, positoque $dQ=Q'dx$, orietur ista aequatio~:}/L
\F{Par cons\'equent, prenant cette diff\'erentielle multipli\'ee par $\tfr x/{dx}$, et posant $dQ=Q'dx$, appara\^\i{}tra cette \'equation ~:}/F
\E{Taking thus this differential multiplied by $\tfr x/{dx}$, and denoting $dQ=Q'dx$, this equation will appear~:}/E
$$
	\fr {x^{n+p}+x^{n-p}}/{(1+x^n)^2} = \fr {Q'x}/{1+x^n}	- \fr {nQx^n}/{(1+x^n)^2} + \fr R/{1+x^n},
 $$
%
%====================================================  (Page 18)
%
\L{quae per $1+x^n$ multiplicata hoc modo repraesentetur~:}/L
\F{qui multipli\'ee par $1+x^n$ se repr\'esente de cette mani\`ere~:}/F
\E{which multiplied by $1+x^n$ can be represented in this fashion~:}/E
$$
	\fr {x^{n+p} + x^{n-p} + nQx^n}/{1+x^n} = Q'x + R,
 $$
\L{ubi iam $Q$ ita accipi debet, ut illa fractio ad integrum revocetur. Facile autem patet, hoc fieri statuendo $n Q = -x^p + x^{n-p}$, tum enim illa fractio fiet}/L
\F{o\`u maintenant $Q$ doit \^etre prise, de sorte que cette expression fractionnaire devienne enti\`ere. Mais il est facile de le faire en posant $n Q = -x^p + x^{n-p}$, car alors cette fraction doit \^etre}/F
\E{where now $Q$ must be taken so that this fractional expression becomes entire. But this is easy to do this by putting $n Q = -x^p + x^{n-p}$, since then this fraction must be}/E
$$
	\fr {x^{n-p}+x^{2n-p}}/{1+x^n} = x^{n-p},
 $$
\L{ita ut nunc habeanus $x^{n-p}=Q'x+R$. Cum igitur sit}/L
\F{de sorte que nous avons $x^{n-p}=Q'x+R$. Comme on a}/F
\E{so that we now have $x^{n-p}=Q'x+R$. Since}/E
$$
	Q = \fr {x^{n-p}-x^p}/n,
 $$
\L{erit}/L
\F{on aura}/F
\E{we will have}/E
$$
	Q' x = \fr {(n-p)x^{n-p} - p x^p}/n,
 $$
\L{hincque colligitur $R=p(x^{n-p}+x^p)/n$, quocirca formula integralis proposita reducta est ad hanc formam~:}/L
\F{par cons\'equent nous obtenons $R=p(x^{n-p}+x^p)/n$, et donc la formule int\'egrale propos\'ee est r\'eduite \`a cette forme~:}/F
\E{thus we get $R=p(x^{n-p}+x^p)/n$, and so the proposed integral formula is reduced to this form~:}/E
$$
	\fr {x^{n-p}-x^p}/{n(1+x^n)} + \fr p/n \, \intofx { x^{n-p}+x^p }/{ 1+x^n }
 $$
\L{quod integrale ita est summendum, ut evanescet posito $x=0$. Nun igitur statuamus {\xon}, ac prior pars absoluta evanescit, formulae autem integralis valor, per ea quae dudum sunt inventa, prodit}/L
\F{cette int\'egrale devant \^etre prise pour s'annuler pour $x=0$. Maintenant nous posons {\xon}, pour que la partie int\'egr\'ee pr\'ec\'edente s'annule, et la valeur de la formule int\'egrale, par ce qui est connu depuis longtemps\footnote{Inst. Calculi Integralis IV -- E660, Supp. III, \S 70, p. 123; ou \S 16 avec $p-n/2,n/2$ -- traducteur.}, devient}/F
\E{this integral being taken so as to vanish for $x=0$. Now we put {\xon}, so that the preceding integrated part vanishes, and the value of the integral formula, from what is known since a long time\footnote{Inst. Calculi Integralis IV -- E660, Supp. III, \S 70, p. 123; or \S 16 with $p-n/2,n/2$ -- translator.}, becomes}/E
$$
	\fr p/n \, \fr \pi/{ n\sin\tfr{p\pi}/n } = \fr { \pi p }/{ nn\sin\tfr{p\pi}/n },
 $$
\L{qui ergo cum ante invento perfecte congruit.}/L
\F{qui s'accorde donc parfaitement avec ce qui a \'et\'e trouv\'e avant.}/F
\E{which thus is in perfect accord with what has been found before.}/E
\medskip
%
%---------------------------------------  (Section VIII)
%
\L{\LP VIII. Tribuantur nunc etiam in hac postrema forma exponenti $p$ valor imaginarius, ponendo $p=q\I$, et cum, ut ante vidimus, hinc fiat $\xpmp=2\cos (q\log x)$, formula integralis proposita erit}/L
\F{\LP VIII. Si l'on attribue dans cette derni\`ere formule des valeurs imaginaires \`a l'exposant $p$, en posant $p=\I q$, et comme, tel que vu pr\'ec\'edemment, ceci fait $\xpmp=2\cos (q\log x)$, la formule int\'egrale propos\'ee sera}/F
\E{\LP VIII. If we attribute in this last formula imaginary values to the exponent $p$, putting $p=\I q$, and since, as seen earlier, this makes $\xpmp=2\cos (q\log x)$, the proposed integral formula will be}/E
$$
	= 2 \intofx { x^n\cos (q\log x) }/{ (1+x^n)^2 }.
 $$
\L{Pro eius valore autem iam ante vidimus fore}/L
\F{Pour sa valeur nous avons d\'ej\`a vu que}/F
\E{For its value we have already seen that}/E
%
%====================================================  (Page 19)
%
$$
	\sin \fr {\I\pi q}/n = \fr { e^{\tfr{-\pi q}/n} - e^{\tfr{+\pi q}/n} }/{ 2\I }
 $$
\L{quamobrem valor nostrae formulae, ab {\xzo} ad \xon, extensus, erit}/L
\F{et ainsi de cette fa\c{c}on la valeur de notre formule, prise de {\xzo} \`a \xon, sera}/F
\E{and thus the value of our formula, taken from {\xzo} to \xon, will be}/E
$$
	\fr { 2\pi q }/{ nn( e^{\tfr{\pi q}/n} - e^{-\tfr{\pi q}/n} ) }
 $$
\L{unde deducimus sequens theorema omni attentione dignum.}/L
\F{et nous d\'eduisons le th\'eor\`eme suivant m\'eritant toute attention.}/F
\E{and we deduce the following theorem worthy of all attention.}/E
\par\medskip\noindent
\L{\cl{ Theorema.}}/L
\F{\cl{ Th\'eor\`eme.}}/F
\E{\cl{ Theorem.}}/E
\par\smallskip\noindent
\L{\cl{ Si valor ipsius formulae integralis~:}}/L
\F{\cl{ La valeur de cette formule int\'egrale~:}}/F
\E{\cl{ The value of the integral formula~:}}/E
$$
	\intofx { x^n\cos (q\log x) }/{ (1+x^n)^2 },
	\fourier{\fr 1/2 \int_0^\infty \fr { \cos(2qt) }/{ \cosh^2(nt) }\,dt}
 $$
\L{\cl{ a termino {\xzo} usque {\xon} extendatur, is semper aequabitur huic formulae:}}/L
\F{\cl{ prise de {\xzo} jusqu'\`a \xon, est toujours \'egale \`a cette formule~:}}/F
\E{\cl{ taken from {\xzo} up to \xon, is always equal to this formula~:}}/E
$$
	\fr { \pi\tfr q/n }/{ n(e^{\pi\tfr q/n} - e^{-\pi\tfr q/n}) }.
	\fourier{\fr 1/{2n} \fr { \pi \tfr q/n }/{ \sinh(\pi \tfr q/n) }}
 $$
\L{Cuius autem Theorematis demonstratio ex principiis iam cognitis vix elici posse videtur.}/L
\F{La d\'emonstration de ce th\'eor\`eme pourrait sembler difficile \`a tirer de principes d\'ej\`a connus.}/F
\E{The demonstration of this theorem could seem difficult to draw from principles already known.}/E
\medskip
%
%---------------------------------------  (Section IX)
%
\L{\LP IX. Praetera etiam perspicuum est, methodum, qua usi sumus ad nostram formulam integrandam, subsistere non posse, nisi terminus medius denominatoris binario sit minor, quam ob caussam eum hac forma $2\cos\theta$ expressimus. Quamobrem hinc oritur quaestio maximi momenti~: utrum nostrae conclusiones etiamnunc valeant, si terminus ille medius binario major acciperetur, sive si angulus $\theta$ foret imaginarius, nec ne ? Verum etiam hoc casu nullum dubium superesse potest, quin formula nostra finalis etiamnunc veritati consentanea sit futura. Ante omnia autem hic est observandum, illi termino medio $2\cos\theta$ valorem negativum tribui convenire, quia alioquin ipse denominator in nihilum abiret, dum quantitas nostra variabilis $x$ a termino $0$ usque ad $1$ augetur. Hanc ob rem statuamus angulum $\theta=\pi-\eta$, et valor noster integralis erit}/L
\F{\LP IX. Par ailleurs il est \'egalement clair que la m\'ethode, que nous avons utilis\'ee pour int\'egrer notre formule, ne peut subsister sans que le terme milieu du d\'enominateur ne soit plus petit, c'est pourquoi nous lui avons donn\'e la forme $2\cos\theta$. Pour cette raison se pose la question importante~: nos conclusions seraient-elles encore valables, si ce terme milieu \'etait pris plus grand, ou si l'angle $\theta$ \'etait imaginaire, ou non ? Mais m\^eme ce cas peut sans doute subsister, mais l'accord final sur la validit\'e de notre formule est \`a venir. Mais avant tout il faut observer ceci, la valeur de ce terme m\'edian $2\cos\theta$ ne peut devenir n\'egative de telle sorte que le d\'enominateur s'annulerait, quand la valeur de notre variable $x$ est augment\'ee de la borne $0$ jusqu'\`a $1$. Pour cette raison posons l'angle $\theta=\pi-\eta$, et la valeur de notre int\'egrale sera}/F
\E{\LP IX. And it is also clear that the method that we have used to integrate our formula cannot function without the middle term being smaller, and that is why we have given it the form $2\cos\theta$. For this reason the important question must be asked~: whether our conclusions would still hold if this middle term would be larger, if the angle $\theta$ was imaginary, or not ? But without doubt even this case can subsist, but the final accord on the validity of our formula is yet to come. Before anything else, however, it must be pointed out that negative values of this middle term $2\cos\theta$ cannot be so attributed that would make the denominator vanish, when the value of our variable $x$ is increased from $0$ to $1$. For this reason let us set the angle $\theta=\pi-\eta$, and the value of our integral will be}/E
%
%====================================================  (Page 20)
%
$$
	\int_{x=0}^{x=1} \fr {dx}/x \fr { \xpmp }/{ x^n+2\cos\eta+x^{-n} }
				= \fr { \pi \sin \tfr{p\eta}/n }/{ n \sin\eta \sin\tfr{p\pi}/n  }.
 $$
\L{In hac igitur formula faciamus angulum $\eta$ imaginarium, ponendo $\eta=\phi\I$, eritque per ea quae iam supra observavimus, $2\cos\phi\I=e^\phi + e^{-\phi}$, ita ut noster denominator sit}/L
\F{Dans cette formule rendons par cons\'equent l'angle $\eta$ imaginaire, en posant $\eta=\phi\I$, et par ce que nous avons observ\'es plus haut, $2\cos\phi\I=e^\phi + e^{-\phi}$, de sorte que notre d\'enominateur soit}/F
\E{Thus in that formula let us make the angle $\eta$ imaginary, defining $\eta=\phi\I$, so that by what we found above, $2\cos\phi\I=e^\phi + e^{-\phi}$, and so that our denominator becomes}/E
$$
	x^n + e^\phi + e^{-\phi} + x^{-n} = \fr 1/{x^n}(x^n+e^\phi)(x^n+e^{-\phi})
 $$
\L{quem idcirco statim in duos factores reales formae $x+k$ resolvere licet; tum vero fiet}/L
\F{qui donc se r\'esout imm\'ediatement en deux facteurs r\'eels de la forme $x+k$; puis effectivement il faut faire}/F
\E{which thus can be resolved  immediately in two real factors of the form $x+k$; then in fact we must make}/E
$$
	\sin \eta = \sin\phi\I = \fr { e^{-\phi} - e^{+\phi} }/{ 2\I },
 $$
\L{similique modo erit}/L
\F{de mani\`ere similaire on aura}/F
\E{similarly we will have}/E
$$
	\sin \fr{p\eta}/n = \sin\fr{p\phi\I}/n 
	= \fr { e^{-\tfr{p\phi}/n} - e^{+\tfr{p\phi}/n} }/{ 2\I },
 $$
\L{unde formula nostra integralis emergit realis}/L
\F{et notre formule int\'egrale devient r\'eelle}/F
\E{and our integral formula becomes real}/E
$$
	= \fr { \pi(e^{-\tfr{p\phi}/n} - e^{+\tfr{p\phi}/n}) }/{
		   n (e^{-\phi} - e^{+\phi}) \sin\tfr{p\pi}/n }.
 $$
\L{Statuamus autem hic brevitatis gratia $e^\phi=f$, ut sit $e^{-\phi}=1/f$, atque nostra formula integralis sequentem induet formam~:}/L
\F{Nous d\'efinissons ici pour abr\'eger $e^\phi=f$, de sorte que $e^{-\phi}=1/f$, et que notre formule int\'egrale acquiert la forme suivante~:}/F
\E{We define here for brevity $e^\phi=f$, so that $e^{-\phi}=1/f$, and our integral formula gets the following form~:}/E
$$
	\int_{x=0}^{x=1} \fr {dx}/x \fr { \xpmp }/{ x^n + f + \tfr 1/f + x^{-n} }
		= \fr { \pi(f^{\tfr p/n}-f^{-\tfr p/n}) }/{ n(f-f^{-1})\sin \tfr{\pi p}/n },
 $$
\L{id quod tanquam theorema omni attentione dignum spectari potest; ubi per se intelligitur, valorem eiusdem integralis, usque ad {\xoo} extensum, fore duplo majorem.}/L
\F{un th\'eor\`eme qui peut \^etre consid\'er\'e comme m\'eritant toute attention; o\`u l'on comprend que la valeur de cette int\'egrale, prise jusqu'\`a {\xoo} serait deux fois plus grande.}/F
\E{a theorem that can be seen worthy of all attention; and it is understood that the value of this integral, taken up to {\xoo}, would be twice as large.}/E
\medskip
%
%---------------------------------------  (Section X)
%
\L{\LP X. Quodsi iam in hac forma etiam exponenti $p$ valorem imaginarium tribuamus, pariter nullo modo dubitari poterit}/L
\F{\LP X. Maintenant si dans cette forme nous donnons \`a l'exposant $p$ des valeurs imaginaires, pareillement on ne peut douter d'aucune fa\c{c}on}/F
\E{\LP X. Now if in this form we give to the exponent $p$ imaginary values, similarly there can be no doubt}/E
%
%====================================================  (Page 21)
%
\L{quin conclusio nostra vera sit mansura. Ponamus igitur $p=q\I$, eritque ut ante $\xpmp=2\cos (q\log x)$; tum vero erit}/L
\F{que notre conclusion demeure valide. Posons dans ce cas $p=q\I$, de sorte que comme avant $\xpmp=2\cos (q\log x)$; alors en effet on aura}/F
\E{that our conclusion remains valid. Let us put in this case $p=q\I$, so that as before $\xpmp=2\cos (q\log x)$; then we will have}/E
$$
	\sin \fr {p\pi}/n = \fr { e^{-\tfr{q\pi}/n} - e^{\tfr{q\pi}/n} }/{ 2\I }
 $$
\L{pro integralis autem numeratore erit}/L
\F{mais pour le num\'erateur de l'int\'egrale on aura}/F
\E{but for the numerator of the integral we will have}/E
$$
	f^{\tfr p/n}-f^{-\tfr p/n} = 2\I\sin \fr {q\log f}/n.
 $$
\L{His igitur valoribus sequens nanciscimur}/L
\F{De ces valeurs on trouve ce qui suit}/F
\E{From these values we find what follows}/E
\par\medskip\noindent
\L{\cl{	Theorema.		}}/L
\F{\cl{	Th\'eor\`eme.	}}/F
\E{\cl{	Theorem.		}}/E
\par\medskip\noindent
\L{\cl{ Valor ipsius formulae integralis~:	}}/L
\F{\cl{ La valeur de cette formule int\'egrale~:}}/F
\E{\cl{ The value of this integral formula~:	}}/E
$$
	\intofx { \cos (q\log x) }/{ x^n + f + \tfr 1/f + x^{-n} },
 $$
\L{\cl{ a termino {\xzo} usque ad {\xon} extensus, semper aequabitur formulae			}}/L
\F{\cl{ prise de la borne {\xzo} jusqu'\`a {\xon}, est toujours \'egale \`a la formule	}}/F
\E{\cl{ taken from {\xzo} up to {\xon}, is always equal to the formula				}}/E
$$
	\fr { 2\pi\sin \tfr {q(\log f)}/n }/{ n(f-f^{-1})(e^{\tfr {q\pi}/n} - e^{-\tfr {q\pi}/n}) }.
 $$
\medskip
%
%---------------------------------------  (Section XI)
%
\L{\LP XI. Deinde iam pridem observaui\footnote{Leonhardi Euleri Inst. Calculi Integralis IV (1845) -- Enestr\"om 660, Supplementum V, \S 188, p. 373.}, omnia huiusmodi integralia satis commode per series infinitas exprimi posse. Cum enim ista fractio~:}/L
\F{\LP XI. Puis, comme nous l'avons observ\'e il y a longtemps\footnote{Leonhard Euler, Inst. Calculi Integralis IV (1845) -- Enestr\"om 660, Supplementum V, \S 188, p. 373.}, toutes les int\'egrales de ce genre peuvent ais\'ement \^etre exprim\'ees par des s\'eries infinies. Comme cette  fraction~:}/F
\E{\LP XI. Then, as we observed long ago\footnote{Leonhard Euler, Inst. Calculi Integralis IV (1845) -- Enestr\"om 660, Supplementum V, \S 188, p. 373.}, all such integrals can be easily expressed by infinite series. Since this  fraction~:}/E
$$
	\fr { x^p }/{ \tortn } = \fr { x^{n+p} }/{ \tornn }
 $$
\L{resolvatur in hanc seriem~:}/L
\F{s'exprime par cette s\'erie~:}/F
\E{can be developed into this series~:}/E
$$
\fr 1/{\sin\theta}\left(x^{n+p}\sin\theta+x^{2n+p}\sin2\theta+x^{3n+p}\sin3\theta + \etc\right)
 $$
\L{integrale istud~:}/L
\F{donc l'int\'egrale~:}/F
\E{thus the integral~:}/E
$$
	\intofx { x^p }/{ \tortn },
	\fourier{\fr 1/2 \int_0^\infty \fr { e^{pt} }/{ \cosh nt - \cos\theta }\,dt}
 $$
\L{a termino {\xzo} usque ad {\xon} extensum, aequabitur huic seriei infinitae~:}/L
\F{prise de la borne {\xzo} jusqu'\`a {\xon}, est \'egale \`a cette s\'erie infinie~:}/F
\E{taken from {\xzo} up to {\xon}, is equal to this infinite series~:}/E
$$
	\fr 1/{\sin\theta}\left(\fr{\sin\theta}/{n+p}+\fr{\sin2\theta}/{2n+p}
		+\fr{\sin3\theta}/{3n+p}+\fr{\sin4\theta}/{4n+p} + \etc\right).
	\fourier{\fr 1/{ n\sin\theta } \sum_{k=1}^\infty \fr { \sin k\theta }/{ k+b }}
 $$
\L{Hinc ergo, si $p$ negative caperemus, tum formula nostra princialis}/L
\F{Ainsi donc, si l'on rend $p$ n\'egatif, alors notre formule principale}/F
\E{And thus, if we make $p$ negative, then our principal formula}/E
$$
	\intofx { \xpmp }/{ \tortn },
	\fourier{\int_0^\infty \fr { \cosh pt }/{ \cosh nt - \cos\theta }\,dt}
 $$
\L{ab {\xzo} ad {\xon} extensa, semper aequabitur huic seriei infinitae geminatae~:}/L
\F{prise de {\xzo} \`a {\xon}, est toujours \'egale \`a cette double s\'erie infinie~:}/F
\E{taken from {\xzo} to {\xon}, is always equal to this double infinite series~:}/E
%
%====================================================  (Page 22)
%
$$
	\fr 1/{\sin\theta}	\left(
	   \fr{\sin\theta}/{n+p} + \fr{\sin2\theta}/{2n+p} + \fr{\sin3\theta}/{3n+p} + \fr{\sin4\theta}/{4n+p}
	 + \etc			\right.
 $$
$$
	 				\left. \qquad\qquad
	 + \fr{\sin\theta}/{n-p} + \fr{\sin2\theta}/{2n-p} + \fr{\sin3\theta}/{3n-p} + \fr{\sin4\theta}/{4n-p}
	 + \etc			\right),
 $$
\L{quae binis homologis conjugendis contrahitur in hanc seriem~:}/L
\F{qui, en groupant les termes semblables deux \`a deux, se contracte dans cette s\'erie~:}/F
\E{which after grouping two by two similar terms contracts into this series~:}/E
$$
	\fr {2n}/{\sin\theta}\left(\fr{\sin\theta}/{nn-pp}+\fr{2\sin2\theta}/{4nn-pp}
		+\fr{3\sin3\theta}/{9nn-pp}+\fr{4\sin4\theta}/{16nn-pp} + \etc\right).
	\hskip -3em
	\fourier{\fr 2/{ n\sin\theta } \sum_{k=1}^\infty \fr { k\sin k\theta }/{ k^2-b^2 }}
 $$
\medskip
%
%---------------------------------------  (Section XIIa)
%
\L{\LP XIIa. Hinc iam manifesto pro casu, quo positur $p=q\I$, ista series infinita exoritur~:}/L
\F{\LP XIIa. Maintenant il est clair qu'en posant $p=q\I$, cette s\'erie infinie s'obtiendra~:}/F
\E{\LP XIIa. Now it is clear that by putting $p=q\I$, this series will be obtained~:}/E
$$
	\fr {2n}/{\sin\theta}\left(\fr { \sin\theta }/{ nn+qq }+\fr { 2\sin2\theta }/{ 4nn+qq }
		+\fr { 3\sin3\theta }/{ 9nn+qq }+\fr { 4\sin4\theta }/{ 16nn+qq } + \etc\right)
	\hskip -3em
	\fourier{\fr 2/{n\sin\theta} \sum_{k=1}^\infty \fr { k\sin k\theta }/{ k^2+(\tfr q/n)^2 }}
 $$
\L{quae ergo exprimit valorem huius formulae integralis~:}/L
\F{qui donc exprime la valeur de cette formule int\'egrale~:}/F
\E{which thus expresses the value of this integral formula~:}/E
$$
	\intofx { 2\cos (q\log x) }/{ \tortn },
	\fourier{\int_0^\infty \fr { \cos qt }/{ \cosh nt - \cos\theta }\,dt}
 $$
\L{scilicet ab {\xzo} ad {\xon} extensae, ita ut istius seriei summa finito modo expressa sit etiam}/L
\F{prise de {\xzo} \`a {\xon}, de sorte que la somme de cette s\'erie s'exprime aussi en termes finis}/F
\E{taken from {\xzo} to {\xon}, so that the sum of this series can be expressed in finite terms}/E
$$
	\fr \pi/{n\sin\theta}\,\left( \fr { 
			e^{-\tfr{q(\pi-\theta)}/n} - e^{+\tfr{q(\pi-\theta)}/n} }/{ 
			e^{-\tfr{q\pi}/n} 	   - e^{+\tfr{q\pi}/n} 		 }\right).
	\fourier{\fr { \pi\sinh(\pi-\theta)\tfr q/n }/{ n\sin(\pi-\theta)\sinh\pi\tfr q/n }}
 $$
\L{Quin etiam facile intelligitur hic quoque angulum $\theta$ imaginarium accipi posse. Vidimus enim posito $\theta=\phi\I$ fore}/L
\F{Par ailleurs il est facile de comprendre que l'angle $\theta$ peut prendre des valeurs imaginaires. Nous avons vu en effet qu'en posant $\theta=\phi\I$ on aurait}/F
\E{And it is easily understood that the angle $\theta$ can take imaginary values. We have seen that after putting $\theta=\phi\I$ we would have}/E
$$
	\sin\theta = \fr { e^{-\phi} - e^{+\phi} }/{ 2\I },
 $$
\L{hincque in genere}/L
\F{et ainsi en g\'en\'eral}/F
\E{and thus in general}/E
$$
	\sin\lambda\theta = \fr { e^{-\lambda\phi} - e^{+\lambda\phi} }/{ 2\I }.
 $$
\L{Quare si statuamus $e^\phi=f$, erit}/L
\F{C'est pourquoi, si nous posons $e^\phi=f$, nous aurons}/F
\E{And so if we put $e^\phi=f$, we will have}/E
$$
	\fr { \sin\lambda\theta }/{ \sin\theta } = 
		\fr { f^\lambda - f^{-\lambda} }/{ f - \tfr 1/f },
 $$
\L{unde series illa satis concinnam formam accipiet.}/L
\F{et cette s\'erie admet une forme assez concise.}/F
\E{and that series has a concise enough form.}/E
\medskip
%
%====================================================  (Page 23)
%
%---------------------------------------  (Section XIIb)
%
\L{\LP XIIb. Denique operationes, quibus in integratione nostrae formulae sumus usi, consistere nequeunt, nisi exponens $n$ fuerit numerus integer. Interim tamen, valor integralis, quem invenimus pro casu vel {\xon} vel {\xoo}, veritati conformis deprehenditur, non solum quando pro $n$ numerus fractus quicunque sed etiam adeo imaginarius accipitur, quorum prius facile ostenditur. Sit enim $n=\tfr m/\lambda$, ac ponatur $x=y^\lambda$, atque ob}/L
\F{\LP XIIb. Enfin les op\'erations, que nous avons utilis\'ees dans nos formules int\'egrales, \'echouent, sauf si l'exposant $n$ est un nombre entier. Cependant, la valeur de l'int\'egrale que nous avons trouv\'ee pour les cas  {\xon} ou {\xoo}, reste valide, non seulement quand $n$ est un nombre fractionnaire quelconque mais aussi quand il prend des valeurs imaginaires, ce qui se montre facilement dans le premier cas. Soit en effet $n=\tfr m/\lambda$, et posons $x=y^\lambda$, de sorte que}/F
\E{\LP XIIb. Finally the operations, which we have used for our integral formulas, fail unless the exponent $n$ is an integer number. However, the value of the integrals that we have found for the cases {\xon} or {\xoo}, remains valid not only when $n$ takes on arbitrary fractional values but also imaginary values, which can be shown easily in the first case. Thus let $n=\tfr m/\lambda$, and put $x=y^\lambda$, so that}/E
$$
	\fr {dx}/x = \lambda \fr {dy}/y,
 $$
\L{orietur haec forma integralis exponentibus integris contenta~:}/L
\F{et appara\^\i{}t cette forme int\'egrale contenant des exposant entiers~:}/F
\E{and there appears this integral form containing integer exponents~:}/E
$$
	\lambda \intofy { y^{\lambda p} + y^{-\lambda p} }/{ y^m-2\cos\theta+y^{-m} }
 $$
\L{cuius ergo valor casu {\xon} debet esse secundum formulam supra inventam}/L
\F{dont la valeur dans le cas {\xon} doit \^etre la deuxi\`eme formule trouv\'ee plus haut}/F
\E{whose value in the case {\xon} must be the second formula found above}/E
$$
	\fr {\lambda\pi}/m\,\fr { \sin \tfr{\lambda p(\pi-\theta)}/m }/{ \sin\theta \sin \tfr{\lambda p\pi}/m },
 $$
\L{qui, si iam loco $m$ valor $\lambda n$ restituatur, manifesto abit in ipsam nostram formulam supra inventam~:}/L
\F{qui, apr\`es remplacement de $m$ par sa valeur $\lambda n$, redonne clairement notre formule trouv\'ee plus haut~:}/F
\E{which, when $m$ is replaced by its value $\lambda n$, clearly gives our formula found above~:}/E
$$
	\fr \pi/n\,\fr { \sin \tfr{p(\pi-\theta)}/n }/{ \sin\theta \sin \tfr{p\pi}/n }.
 $$
\L{Hinc autem nulli amplius dubio relinquitur, quin veritas haec subsistat, etiamsi $n$ fuerit numerus imaginarius\footnote{Manifestum est hanc conclusionem falsam esse, quotiescumque pro $n$ numerius  imaginarius formae $m\I$ accipitur. (Alexander Liapounoff, 1920)}. Ponamus igitur $n=m\I$; et formula integralis reducetur ad hanc formam~:}/L
\F{D'autre part il ne reste aucun doute, que leur validit\'e demeure, quand $n$ prend des valeurs imaginaires\footnote{Il est manifeste que cette conclusion est fausse, d\`es que l'on prend pour $n$ un nombre imaginaire de la forme $\I m$. (Alexander Liapounoff, 1920)}. Posons donc $n=m\I$; et la formule int\'egrale sera r\'eduite \`a cette forme~:}/F
\E{In the other case there is no doubt that their validity remains, when $n$ takes imaginary values\footnote{It is manifest that this conclusion is false, whenever we take for $n$ an imaginary number of the form $\I m$. (Alexander Liapounoff, 1920)}. Thus, let us put $n=m\I$; and the integral formula will be reduced to this form~:}/E
$$
	\intofx { \xpmp }/{ 2\cos (m\log x) - 2 \cos\theta },
 $$
\L{cuius ergo valor casu {\xon} certe erit}/L
\F{dont la valeur pour le cas {\xon} sera certainement}/F
\E{whose value for the case {\xon} will certainly be}/E
$$
\fr \pi/{m\I}\, \fr { e^{\tfr{p(\pi-\theta)}/m} - e^{-\tfr{p(\pi-\theta)}/m} }/{
		 \sin\theta(e^{\tfr{p\pi}/m} - e^{-\tfr{p\pi}/m}) }
 $$
\L{ubi mirum videbitur istum valorem semper esse imaginarium, licet ipsa formula differentialis, dum variabilis $x$, a termino $0$ usque ad terminum {\xon} augetur, maneat realis, id quod merito maxime videtur paradoxum. Interim tamen non desunt casus, quibus valor integralis formulae differentialis realis manifesto evadit imaginarius, id quod in ista formula simpliciori}/L
\F{o\`u il sera consid\'er\'e \'etonnant que cette valeur soit toujours imaginaire, bien que cette forme diff\'erentielle, pour la variable $x$, de la borne {\xzo} jusqu'\`a la borne {\xon}, reste r\'eelle, ce qui m\'erite d'\^etre vu comme un grand paradoxe. Dans l'intervalle cependant ne manquent pas les cas dans lesquels la valeur de l'int\'egrale de diff\'erentielles r\'eelles devient manifestement imaginaire, ce qui dans cette formule simple}/F
\E{where it will be considered surprising that this value be always imaginary, although this differential form, for the variable $x$, from the limit {\xzo} up to the limit {\xon}, stays real, which merits to be seen as a great paradox. Meanwhile however cases are not lacking where the value of an integral of real differentials becomes clearly imaginary, which in this simple formula}/E
$$
	\int \fr { dx }/{ x\cos (m\log x) }
 $$
\medskip
%
%====================================================  (Page 24)
%
\L{ostendisse sufficiet, quae utique, dum $x$ a $0$ ad $1$ augetur, constanter manet realis. Ad hanc ergo formulam integrandam statuamus $\log x=-z$, ubi notetur, dum $x$ a $0$ usque ad $1$ progreditur, tum quantitatem $z$ ab $\infty$ usque ad $0$ decrescere. Nunc igitur formula nostra integralis erit}/L
\F{est d\'emontr\'e suffisamment, quoique bien s\^ur, quand  $x$ augmente de $0$ \`a $1$, elle demeure constamment r\'eelle. Pour int\'egrer cette formule posons $\log x=-z$, o\`u il faut noter, quand $x$ progresse de $0$ jusqu'\`a $1$, alors la quantit\'e $z$ d\'ecroit de $\infty$ jusqu'\`a $0$. Maintenant notre formule int\'egrale sera}/F
\E{is sufficiently demonstrated, although of course, when $x$ increases from $0$ to $1$, it stays constantly real. To integrate this formula, let $\log x=-z$, where it must be noted that, when $x$ progresses from $0$ up to $1$, then the quantity $z$ decreases from $\infty$ to $0$. Now our integral formula will be}/E
$$
	\int \fr { - dz }/{ \cos mz },
 $$
\L{cum vero constet esse}/L
\F{comme on constate que}/F
\E{as we see that}/E
$$
	\int \fr { d\phi }/{ \sin\phi } = \log \tan \fr \phi/2,
 $$
\L{sumamus $\phi=90^\circ-mz$, eritque $d\phi=-mdz$, hincque}/L
\F{prenons $\phi=90^\circ-mz$, on aura aussi $d\phi=-mdz$, et}/F
\E{let us take $\phi=90^\circ-mz$, we will have $d\phi=-mdz$, and}/E
$$
	\int \fr { - m dz }/{ \cos mz } = + \log \tan(45^0-\fr {mz}/2),
 $$
\L{quod integrale manifesto evanescit pro termino $z=0$, dum autem ab hoc termino quantitas $z$ in infinitum usque augetur, infinites tangens huius anguli fiet negativa, eiusque logarithmus propterea imaginarius, unde non amplius mirabimur, quod formulae differentialis realis integrale evadere possit certis casibus imaginarium.}/L
\F{cette int\'egrale s'annulant manifestement pour la borne $z=0$, mais la quantit\'e $z$ croissant \`a partir de cette borne jusqu'\`a l'infini, la tangente de l'angle \'etant n\'egative, et le logarithme devenant imaginaire, donc il n'est plus \`a se demander comment des int\'egrales de formules diff\'erentielles r\'eelles deviennent dans certains cas imaginaires.}/F
\E{this integral vanishing for the limit $z=0$, but the quantity $z$ increasing from this limit to infinity, the tangent of the angle being negative and the logarithm being imaginary, so we do not need to ask the question of how integrals of real differential formulas become imaginary in certain cases.}/E
\medskip
%
%---------------------------------------  (Section XIII)
%
\L{\LP XIII. Hoc igitur modo evictum est formulae nostrae differentialis propositae}/L
\F{\LP XIII. De cette fa\c{c}on nous sommes assur\'es que notre formule diff\'erentielle propos\'ee}/F
\E{\LP XIII. In this manner, we are assured that our proposed differential formula}/E
$$
	\intofx { \xpmp }/{ \tortn },
	\fourier{= \fr { \pi\sin\tfr {(\pi-\theta)p}/n }/{ n\sin(\pi-\theta)\sin\tfr {\pi p}/n }}
 $$
\L{integrale assignatum a termino {\xzo} usque ad {\xon} semper cum veritate consistere, quicunque valores ternis litteris $n, p$ et $\theta$, tribuantur, sive integri, sive fracti, sive etiam imaginarii. Interim tamen dantur casus iam initio indicati, quibus isti valores integrales a veritae aberrabunt, quippe quod semper usu venire debet, quoties exponens $p$ major est exponente $n$, quam ob causam sedulo excludere debemus omnes casus, quibus formula $p-n$ evadit realis et positiva. His autem exceptis variae formulae, ad quas hic sumus perducti, ita sunt comparatae, ut maxima attentione dignae videantur, simulque non contemnenda incrementa Scientae analyticae promittant.}/L
\F{pour l'int\'egrale prise de {\xzo} jusqu'\`a {\xon} est toujours valide, quelles que soient les valeurs attribu\'ees aux trois lettres $n, p$ et $\theta$, soit enti\`eres, soit fractionnaires, soit m\^eme imaginaires. En m\^eme temps, cependant, dans les cas qui sont indiqu\'es plus haut, sont exclues les valeurs qui rendent aberrante la valeur de l'int\'egrale chaque fois que l'exposant $p$ est plus grand que l'exposant $n$ et nous devons aussi exclure tous les cas o\`u la partie r\'eelle de $p-n$ est positive. Mais cela except\'e, nous sommes arriv\'es \`a des formules diverses qui peuvent \^etre consid\'er\'ees comme m\'eritant la plus grande attention, sans n\'egliger dans le m\^eme temps l'augmentation des connaissances pour la science de l'analyse.}/F
\E{for the integral taken from {\xzo} to {\xon} is always valid, for any values attributed to the three letters $n, p$ and $\theta$, either integers, fractions, or even imaginary numbers. At the same time, however, in the cases which have been indicated above, are excluded the values which render aberrant the value of the integral each time that the exponent $p$ is larger than the exponent $n$ and we must also exclude all cases where the real part of $p-n$ is positive. But that excepted, we have found diverse formulas that can be considered to be worthy of the greatest attention, without neglecting the increase in knowledge for the science of analysis.}/E
%
%	Force inclusion of list of changes in Latin version
%
%Latin: \def\F#1/F{\par\noindent #1}
%Latin: \def\E#1/E{\par\noindent #1}
%
%=================================================================================================
%
\vfill
\eject
\F{\cl{	Annexe A -- Modifications pour la traduction du latin			}}/F
\E{\cl{	Appendix A -- Modifications for the translation from the Latin	}}/E
\bigskip\noindent
\F{1. $\I$ remplace $\sqrt{-1}$ (\TeX\ macro I).}/F
\E{1. $\I$ replaces $\sqrt{-1}$ (\TeX\ macro I).}/E
\par\bigskip\noindent
\F{2. $k$ remplace $i$ lorsque utilis\'e comme indice.}/F
\E{2. $k$ replaces $i$ when used as an index.}/E
\par\bigskip\noindent
\F{3. $\log$ remplace $l$ pour le logarithme naturel.}/F
\E{3. $\log$ replaces $l$ for the natural logarithm.}/E
\par\bigskip\noindent
\F{4. $dx$ remplace $\partial x$ comme diff\'erentielle.}/F
\E{4. $dx$ replaces $\partial x$ as differential.}/E
\par\bigskip\noindent
\F{5. $/4\sin^2(.)$ remplace $/4\sin(.)^2$ (\S 12, \S 13, \S 14).}/F
\E{5. $/4\sin^2(.)$ replaces $/4\sin(.)^2$ (\S 12, \S 13, \S 14).}/E
\par\bigskip\noindent
\F{6. $-b\cos(\alpha+2n\beta)$ remplace $+b\cos(\alpha+2n\beta)$ (derni\`ere \'equation, \S 13).}/F
\E{6. $-b\cos(\alpha+2n\beta)$ replaces $+b\cos(\alpha+2n\beta)$ (last equation, \S 13).}/E
\par\bigskip\noindent
\F{7. $Q=V=(-.+.)/.$ remplace $Q=V=-(.+.)/.$ (avant derni\`ere \'equation, \S 14).}/F
\E{7. $Q=V=(-.+.)/.$ replaces $Q=V=-(.+.)/.$ (second to last equation, \S 14).}/E
\par\bigskip\noindent
\F{8. $x^n+2+x^{-n}$ remplace $x^n-1+x^{-n}$ (th\'eor\`eme, \S VII). }/F
\E{8. $x^n+2+x^{-n}$ replaces $x^n-1+x^{-n}$ (theorem, \S VII).}/E
\par\bigskip\noindent
\F{9. $dQ=Q'dx$ remplace $dQ'=Qdx$ (d\'emonstration, \S VII).}/F
\E{9. $dQ=Q'dx$ replaces $dQ'=Qdx$ (demonstration, \S VII).}/E
\par\bigskip\noindent
\F{10. $x^{\I q}+x^{-\I q}$ remplace $x^p+x^{-p}$ (th\'eor\`eme, \S VIII).}/F
\E{10. $x^{\I q}+x^{-\I q}$ replaces $x^p+x^{-p}$ (theorem, \S VIII).}/E
\par\bigskip\noindent
\F{11. $\tortn$ remplace $\tornn$ (\S XI).}/F
\E{11. $\tortn$ replaces $\tornn$ (\S XI).}/E
\par\bigskip\noindent
\F{12. $\sin\eta$ remplace $\sin\theta$ (\S XI).}/F
\E{12. $\sin\eta$ replaces $\sin\theta$ (\S XI).}/E
\par\bigskip\noindent
\F{13. $2\cos (q\log x)$ remplace $\cos (q\log x)$ (\S XIIa).}/F
\E{13. $2\cos (q\log x)$ replaces $\cos (q\log x)$ (\S XIIa).}/E
\par\bigskip\noindent
\F{14. $\tfr \pi/{m\I}$ remplace $\tfr p/{m\I}$ (\S XIIb).}/F
\E{14. $\tfr \pi/{m\I}$ replaces $\tfr p/{m\I}$ (\S XIIb).}/E
\par\bigskip\noindent
\F{Les corrections 5--13 ont aussi \'et\'e faites par Alexander Liapounoff lors de la r\'e\'edition en 1920 de l'article dans le volume XVIII de la s\'erie 1 des Opera Omnia. L'article en latin y est pr\'efac\'e du sommaire de 1788 en fran\c{c}ais.}/F
\E{The corrections 5--13 have also been made by Alexander Liapounoff for the re-edition in 1920 of the article in volume XVIII of the series 1 of the Opera Omnia. The article in Latin is prefaced there by the summary of 1788 in French.}/E
%
%	End of Latin version
%Latin: \Latex{\end{document}}
%Latin: \bye
%
%=================================================================================================
%
\vfill
\eject
\F{\cl{	Annexe B -- Sommaire et preuves par Sim\'eon Denis Poisson [a]		}}/F
\E{\cl{	Appendix B -- Summary and proofs by Sim\'eon Denis Poisson [a]		}}/E
\bigskip
\F{\LP 22. Les formules les plus simples qui soient comprises dans les pr\'ec\'edentes, sont celles que l'on obtient en faisant $\theta=0$ dans les \'equations (3) et (7). Si l'on met en m\^eme temps $2m$ \`a la place de $m$ dans la premi\`ere, et que l'on fasse les r\'eductions, on a\footnote{Poisson a utilis\'e les variables $p, m, \theta$ au lieu de $m, p, a$ ici, et r\'ef\`ere aussi \`a Plana -- traducteur.}}/F
\E{\LP 22. The simplest formulas included in the preceding ones are those obtained in putting $\theta=0$ in the equations (3) and (7). If at the same time $2m$ is substituted instead of $m$ in the first one, and that the reductions are done, one obtains\footnote{Poisson used the variables $p, m, \theta$ instead of $m, p, a$ here, and refers also to Plana -- translator.}}/E
$$
	\int \fr {e^{2mt}-e^{-2mt}}/{e^{\pi t}-e^{-\pi t}} \, dt = \fr 1/2 \, \tan m  ;
 $$
$$
	\int \fr {e^{2mt}+e^{-2mt}}/{e^{\pi t}+e^{-\pi t}} \, dt = \fr 1/2 \, \sec m .
 $$
\F{Si l'on d\'esigne par $x$ une nouvelle variable; par $n$ un exposant positif; qu'on fasse $e^{-\pi t}=x^n$, et, pour abr\'eger, $\tfr {2mn}/\pi=p$; les valeurs extr\^emes $t=0$ et $t=\tfr 1/0$, r\'epondront \`a {\xon} et \`a {\xzo}, et nos deux \'equations deviendront}/F
\E{If we denote by $x$ a new variable; by $n$ a positive exponent; we substitute $e^{-\pi t}=x^n$, and, for brevity, $\tfr {2mn}/\pi=p$; the extreme values $t=0$ and $t=\tfr 1/0$, will correspond to {\xon} and to {\xzo}, and our two equations will become}/E
$$
	\int \fr {x^{p}-x^{-p}}/{x^n-x^{-n}} \, \fr {dx}/x = \fr\pi/{2n} \, \tan\fr {\pi p}/{2n}  ;
 $$
$$
	\int \fr {x^{p}+x^{-p}}/{x^n+x^{-n}} \, \fr {dx}/x = \fr\pi/{2n} \, \sec\fr {\pi p}/{2n}  ;
 $$
\F{les int\'egrales \'etant prises depuis {\xzo} jusqu'\`a {\xon}.}/F
\E{the integrals being taken from {\xzo} to {\xon}.}/E
\F{Ces formules sont dues, comme on sait, \`a Euler, qui les a d\'emontr\'ees pour le cas seulement o\`u l'exposant $p$ est r\'eel; l'analyse pr\'ec\'edente prouve qu'elles ont encore lieu, lorsque cet exposant est suppos\'e imaginaire.}/F
\E{These formulas are due, it is known, to Euler, who proved them only for the case when the exponent $p$ is real; the preceding analysis proves that they are still valid when that exponent is supposed imaginary.}/E
\bigskip
\F{\LP 28. Par le proc\'ed\'e direct de l'int\'egration des fractions rationnelles, Euler est parvenu \`a ce r\'esultat remarquable par sa simplicit\'e:}/F
\E{\LP 28. By the direct process of integration of rational functions, Euler obtained this result, noteworthy by its simplicity:}/E
$$
	\int \fr {x^{p}+x^{-p}}/{x^n+2\cos a+x^{-n}} \, \fr {dx}/x = 
	\fr { \pi\sin\tfr {a p}/n }/{ n\sin a\sin\tfr {\pi p}/n } ;
 $$
\F{l'int\'egrale \'etant prise depuis {\xzo} jusqu'\`a {\xon}, $a$ \'etant plus petit que $\pi$, et $p$ et $n$ d\'esignant des nombres entiers et positifs, dont le premier est moindre que le second. Il observe ensuite que la m\^eme formule subsiste, lorsque ces exposans sont des quantit\'es r\'eelles quelconques; car si l'on met $x^q$ \`a la place de $x$, $q$ \'etant une quantit\'e r\'eelle et positive, on ne change rien aux limites de l'int\'egrale, et en faisant $qn=n'$, $qp=p'$, il vient}/F
\E{the integral being taken from {\xzo} to {\xon}, $a$ being less than $\pi$, and $p$ and $n$ denoting positive integers such that the first is less than the second. He then observes that the same formula subsists when these exponents are arbitrary real quantities; since if $x^q$ is substituted for $x$, $q$ being a real and positive quantity, nothing is changed for the limits of integration, and defining $qn=n'$, $qp=p'$, we obtain}/E
$$
	\int \fr {x^{p'}+x^{-p'}}/{x^{n'}+2\cos a+x^{-n'}} \, \fr {dx}/x = 
	\fr { \pi\sin\tfr {a p'}/n' }/{ n'\sin a\sin\tfr {\pi p'}/n' } ;
 $$
\F{et \`a cause du facteur ind\'etermin\'e $q$, on peut prendre maintenant pour $p'$ et $n'$ des quantit\'es r\'eelles et positives quelconques. Euler est aussi conduit, par l'induction fond\'ee sur la g\'en\'eralit\'e des formules analytiques, \`a regarder l'exposant $p$ comme imaginaire; en le rempla\c{c}ant donc par $\I q$, il obtient cet autre r\'esultat (M\'emoires de Petersbourg, ann\'ee 1785 et 1787) :}/F
\E{and because of the indeterminate factor $q$, one can now take for $p'$ and $n'$ arbitrary real positive quantities. Euler is also lead, by induction based on the generality of analytic formulas, to view the exponent $p$ as imaginary; replacing it by $\I q$, he deduces this other result (Petersbourg Memoirs, years 1785 and 1787) :}/E
$$
	\int \fr {\cos(q\log x)}/{x^n+2\cos a+x^{-n}} \, \fr {dx}/x = 
	\fr \pi/{ 2 n \sin a }\,
	\fr{ e^{\fr {a q}/n } - e^{- \fr {a q}/n } }/
		{ e^{\fr {\pi q}/n } - e^{- \fr {\pi q}/n } },
%	\fr { \pi\sinh\tfr {a q}/n }/{ n\sin a\sinh\tfr {\pi q}/n } ;
 $$
\F{qu'il ne donne que comme une simple conjecture, en observant qu'il serait \`a d\'esirer qu'on p\^ut trouver un moyen direct d'y parvenir et d'en constater l'exactitude. Or, si l'on fait, dans l'\'equation (15), $e^{-t}=x^n$ et $kn=q$, les limites de l'int\'egrale relative \`a $x$ seront {\xzo} et {\xon}, et cette \'equation se changera dans la pr\'ec\'edente, qui se trouvera, de cette mani\`ere, rigoureusement d\'emontr\'ee.}/F
\E{which he gives as a simple conjecture, observing that it would be desirable that someone would find a direct method to obtain it and verify its validity. But if we substitute, in equation (15), $e^{-t}=x^n$ and $kn=p$, the limits of integration relative to $x$ become {\xzo} and {\xon}, and this equation will be identical to the preceding one, which will, in this manner, be rigorously demonstrated.}/E
\par\bigskip\noindent
\F{[a] S.D. Poisson, Journal de l'\'Ecole Polytechnique, Cahier 18/Tome XI (1820), \S 28, p. 298, 308.}/F
\E{[a] S.D. Poisson, Journal de l'\'Ecole Polytechnique, Cahier 18/Tome XI (1820), \S 28, p. 298, 308.}/E
%
%=================================================================================================
%
\vfill
\eject
\F{\cl{	Annexe C -- Court sommaire par Heinrich Burkhardt [b]		}}/F
\E{\cl{	Appendix C -- Short summary by Heinrich Burkhardt [b]		}}/E
\bigskip
\F{L'\'equation}/F
\E{The equation}/E
$$
	\int_0^\infty \fr{\cos bt \, dt}/{\cosh t + \cos a} = \fr \pi/{\sin a} \fr{\sinh ab}/{\sinh \pi b}
 $$
\F{a d\'ej\`a \'et\'e prouv\'ee par Euler [c] en utilisant une m\'ethode diff\'erente. Il obtient tout d'abord pour des valeurs enti\`eres de $p$ et $n$, en effectuant l'int\'egration ind\'efinie et quelques r\'eductions, l'\'equation\footnote{Burkhardt  a utilis\'e les variables $x, m, r$ au lieu de $b, p, x$ ici, et a oubli\'e le facteur 2 \`a droite -- traducteur.}}/F
\E{has already been proven by Euler [c] using a different method. He first obtains for integer values of $p$ and $n$, by performing indefinite integration and some reductions, the equation\footnote{Burkhardt used the variables $m,r$ instead of $p,x$ used here, and forgot the factor 2 on the right hand side -- translator.}}/E
$$
	\int_0^\infty \fr {x^{p}+x^{-p}}/{x^{n} + 2\cos a + x^{-n}} \, \fr {dx}/x = 
			  \fr {2\,\pi}/{n\sin a} \fr{\sin a p/n}/{\sin \pi p/n} \, ;
	\hskip -3em
	\fourier{2\int_0^1 \fr {x^{p}+x^{-p}}/{x^{n} + 2\cos a + x^{-n}}\fr {dx}/x}
 $$
\F{il fait ensuite remarquer qu'elle reste inchang\'ee quand on prend pour $p$ et $n$ des valeurs fractionnaires partageant le m\^eme d\'enominateur entier\footnote{En fait Euler prouve plus, r\'eduisant le cas $p$ r\'eel et $n$ entier au cas $p$ entier et $n$ entier (\S III); puis le cas $p$ r\'eel et $n$ r\'eel au cas $p$ r\'eel et $n$ entier (\S XIIb), si nous lui pr\^etons la justification requise pour permuter limite et int\'egration sur $[0,1]$, tel que demand\'e dans sa derni\`ere phrase du \S III -- traducteur.}, et il soutient qu'elle est valide pour des valeurs irrationnelles et m\^eme imaginaires de $p$. On obtient la premi\`ere \'equation plus haut par la substitution $x^n=e^{t}$.}/F
\E{then he points out that it remains unchanged when $p$ and $n$ are given fractional values sharing the same integer denominator\footnote{In fact Euler proves more, reducing the case $p$ real and $n$ integer to the case $p$ integer and $n$ integer (\S III); then the case $p$ real and $n$ real to the case $p$ real and $n$ integer (\S XIIb), if we lend him the justification required to interchange limit and integration over $[0,1]$, as requested in his last sentence of \S III -- translator.}, and he claims that it is valid for irrational and even for imaginary values of $p$. The first equation above is obtained after the substitution $x^n=e^{t}$.}/E
\par\bigskip\noindent
\F{[b] H. Burkhardt, Trigonomische Reihen und Integrale (bis etwa 1850), Encyclop\"aedie der Mathematischen Wissenschaften, Zweiter Band, Erster Teil, Zweite H\"alfte, 1904-1916, p. 1133. (Cite aussi Poisson et Plana)}/F
\E{[b] H. Burkhardt, Trigonomische Reihen und Integrale (bis etwa 1850), Encyclop\"aedie der Mathematischen Wissenschaften, Zweiter Band, Erster Teil, Zweite H\"alfte, 1904-1916, p. 1133. (Quotes also Poisson and Plana)}/E
\par\bigskip\noindent
\F{[c] Petrop. n. a. 3 (1785[88]), p. 14 (de 1776); r\'esultat sans preuve dans 5 (1785/89) p. 14 (aussi de 1776). Euler d\'eduit en outre de ce th\'eor\`eme int\'egral la d\'ecomposition en fractions partielles de certaines fonctions.}/F
\E{[c] Petrop. n. a. 3 (1785[88]), p. 14 (from 1776); result without proof in 5 (1785/89) p. 14 (also from 1776). Euler extracts furthermore from this integral theorem the partial fraction decomposition of certain functions.}/E
%
%=================================================================================================
%
\vfill
\eject
\F{\cl{	Annexe D -- Liste de formules d'int\'egration	}}/F
\E{\cl{	Appendix D -- List of integration formulas	}}/E
\PN     %--------------------------------------
\F{Certaines formules int\'egrales de l'article d'Euler sont pr\'esent\'ees ici en utilisant la notation de Fourier avec bornes d'int\'egration et en transformant l'intervalle $(0,1]$ \`a $[0,\infty)$ par la substitution de Poisson $x^n=e^{-t}$.  Suivant Poisson et Burkhardt nous d\'efinissons aussi, pour simplifier, $a=\pi-\theta$, $c=\pi-\zeta$ et $b=\tfr p/n$.}/F
\E{Some integral formulas found in Euler's article are presented here using the Fourier notation with limits of integration,  after transforming the interval $(0,1]$ to $[0,\infty)$ with the Poisson substitution $x^n=e^{-t}$. Following Poisson and Burkhardt we also define, in order to simplify, $a=\pi-\theta$, $c=\pi-\zeta$ and $b=\tfr p/n$.}/E
\PN     %-------------------------------------- 1-15
\F{\S 1--\S 15. Formule principale\footnote{Bierens de Haan (1867) 6.20; Gradshteyn\&Ryzhik (2007) 3.514.2, $0<a<\pi$ et $0<|b|<1$.} valide, par continuit\'e si $a=0$ ou $b=0$, pour $|\Re{a}|<\pi$ et $|\Re{b}|<1$~:}/F
\E{\S 1--\S 15. Main formula\footnote{Bierens de Haan (1867) 6.20; Gradshteyn\&Ryzhik (2007) 3.514.2, $0<a<\pi$ and $0<|b|<1$.} valid, by continuity if $a=0$ or $b=0$, for $|\Re{a}|<\pi$ and $|\Re{b}|<1$~:}/E
$$
	S(a,b,c) :=
	\int_0^\infty \fr { \cosh bt + \cos c }/{ \cosh t + \cos a } \, dt =
		\fr 1/2 \int_{-\infty}^\infty \fr { \cosh bt + \cos c }/{ \cosh t + \cos a } \, dt =
		\fr { \pi\sin ab }/{ \sin a\sin \pi b } + \fr { a\cos c }/{ \sin a }.
 $$
\PN     %-------------------------------------- 6
\F{\S 6. La preuve d'Euler pour for $-\pi<a<\pi$ et $-1<b<1$ utilise les fractions partielles conduisant \`a la formule \'el\'ementaire d'int\'egration ind\'efinie suivante qui se v\'erifie par diff\'erentiation,}/F
\E{\S 6. The proof by Euler for $-\pi<a<\pi$ and $-1<b<1$ uses partial fractions leading to the following elementary indefinite integration formula which can be verified by differentiation,}/E
$$
	\int_0^x \fr { \sin a }/{ y^2 + 2y\cos a + 1 } \, dy =
	\atan \left( \fr { x\sin a }/{ 1 + x\cos a } \right),
 $$
\F{et qui donne cette int\'egrale d\'efinie\footnote{Bierens de Haan 27.22; Gradshteyn\&Ryzhik 3.514.1; Brychkov\&Marichev\&Prudnikov (1986) 2.3.14.29.}${}^,$\footnote{Suivant Weierstrass, cette \'equation, prouv\'ee deux fois par Euler pour $-\pi<a<\pi$ (\S 9, \S I), reste valide pour $|\Re{a}|<\pi$ puisque l'int\'egrale par convergence uniforme et $\tfr a/{\sin a}$ y sont tous deux analytiques; leur diff\'erence est nulle dans toute r\'egion de r\'egularit\'e intersectant $-\pi<a<\pi$ (E.C. Titchmarsh, The Theory of Functions, 2nd ed., 1939, \S 4.4: Analytic Continuation. Integrals containing a complex parameter, p. 147).} en substituant $y=e^{-t}$ et $x=1$ (premi\`ere int\'egrale du \S I)~:}/F
\E{and which gives this definite integral\footnote{Bierens de Haan 27.22; Gradshteyn\&Ryzhik 3.514.1; Brychkov\&Marichev\&Prudnikov (1986) 2.3.14.29.}${}^,$\footnote{Following Weierstrass, this equation, proven twice by Euler for $-\pi<a<\pi$ (\S 9, \S I), remains valid for $|\Re{a}|<\pi$ since the integral by uniform convergence and $\tfr a/{\sin a}$ are analytic there; their difference vanishes in any region of regularity intersecting $-\pi<a<\pi$ (Edward Charles Titchmarsh, The Theory of Functions, 2nd ed., Oxford, 1939, \S 4.4: Analytic continuation. Integrals containing a complex parameter, p. 147).} by substituting $x=1$ and $y=e^{-t}$ (first integral of \S I)~:}/E
$$
	S(a,0,\fr \pi/2)\,=\,\int_0^\infty \fr 1/{ \cosh t + \cos a }\,dt = \fr a/{ \sin a },\qquad |\Re{a}|<\pi.
 $$
\PN     %-------------------------------------- 16
\F{\S 16. Un second cas particulier de la formule principale, $a=c=\tfr \pi/2$, est le plus simple~:}/F
\E{\S 16. A second particular case of the main formula, $a=c=\tfr \pi/2$, is the simplest~:}/E
$$
	S(\fr \pi/2,b,\fr \pi/2) =
	\int_0^\infty \fr { \cosh bt }/{ \cosh t }\,dt = \fr \pi/2 \sec\fr{\pi b}/2 ,\qquad |\Re{b}|<1.
 $$
\F{Un troisi\`eme cas particulier\footnote{dont la validit\'e, pour tout $|\Re{a}|<\pi$ fix\'e, se prolonge encore de $-1<b<1$ \`a $|\Re b|<1$.}, $c=\tfr \pi/2$, est \'etudi\'e dans le reste de l'article sous le symbole $P$~:}/F
\E{A third particular case\footnote{whose validity, for each fixed $|\Re{a}|<\pi$, can be extended again from $-1<b<1$ to $|\Re b|<1$.}, $c=\tfr \pi/2$, is studied in the rest of the article under the symbol $P$~:}/E
$$
	P(a,b) := S(a,b,\fr \pi/2) =
	\int_0^\infty \fr { \cosh bt }/{ \cosh t + \cos a }\,dt = \fr { \pi\sin ab }/{ \sin a\sin \pi b },
	\qquad |\Re{a}|<\pi,\, |\Re{b}|<1.
 $$
\PN     %-------------------------------------- IV
\F{\S IV. Une transform\'ee cosinus de Fourier (de la s\'ecante hyperbolique si $a=\tfr\pi/2$) est obtenue de la formule $P$ en rempla\c{c}ant formellement $b$ par $\I b$~:}/F
\E{\S IV. A Fourier cosine transform (of the hyperbolic secant if $a=\tfr\pi/2$) is obtained from formula $P$ by replacing formally $b$ with $\I b$~:}/E
$$
	P(a,\I b) =
	\int_0^\infty \fr { \cos bt }/{ \cosh t + \cos a }\,dt = \fr { \pi\sinh ab }/{ \sin a \sinh \pi b},
	\qquad |\Re{a}|<\pi,\, |\Im{b}|<1.
 $$
\PN     %-------------------------------------- V
\F{\S V. Une transform\'ee de Laplace bilat\'erale (de la s\'ecante hyperbolique si $a=\tfr\pi/2$) suit par parit\'e~:}/F
\E{\S V. A two sided Laplace transform (of the hyperbolic secant if $a=\tfr\pi/2$) follows using parity~:}/E
$$
	\int_{-\infty}^\infty \fr {      e^{\pm bt}     }/{ \cosh t + \cos a }\,dt =
	\int_{-\infty}^\infty \fr { \cosh bt \pm \sinh bt }/{ \cosh t + \cos a }\,dt =
	\fr { 2 \pi\sin ab }/{ \sin a \sin \pi b},
	\quad |\Re{a}|<\pi,\, |\Re{b}|<1.
 $$
\PN     %-------------------------------------- VI
\F{\S VI. Un premier paradoxe porte sur la p\'eriodicit\'e en $a$ de $\sin ab$, ignorant la condition $|\Re{a}|<\pi$~:}/F
\E{\S VI. A first paradox concerns the periodicity in $a$ of $\sin ab$, ignoring the condition $|\Re{a}|<\pi$~:}/E
$$
	P(a,b) =
	\int_0^\infty \fr { \cosh bt }/{ \cosh t + \cos(a\pm 2k\pi) }\,dt =
	P(a\pm 2k\pi,b) =	\fr{ \pi\sin(a\pm 2k\pi)b }/{ \sin a\sin \pi b } \quad ??
  $$
\PN     %-------------------------------------- VII
\F{\S VII. Cas limite $a\to0$ de $P$, que Euler v\'erifie par r\'eduction au premier cas particulier~:}/F
\E{\S VII. Limit case $a\to0$ of formula $P$, which Euler verifies by reducing it to the first particular case~:}/E
$$
	P(0,b) = \int_0^\infty \fr { \cosh 2bt }/{ \cosh^2 t }dt = 
	\left. \pm \fr {\sinh(2b \pm 1)t}/{\cosh t}\right|_0^\infty	\mp 2 b\, P(\fr\pi/2,2b \pm 1) =
	\fr { \pi b }/{ \sin \pi b },
	\quad |\Re{b}|<1.
  $$
\PN     %-------------------------------------- VIII
\F{\S VIII. La transform\'ee cosinus de Fourier de la s\'ecante hyperbolique au carr\'e est obtenue en rempla\c{c}ant formellement $b$ par $\I b$ dans la formule pr\'ec\'edente~:}/F
\E{\S VIII. The Fourier cosine transform of the hyperbolic secant squared is obtained by replacing formally $b$ with $\I b$ in the last formula~:}/E
$$
	P(0,\I b) =
	\int_0^\infty \fr { \cos 2bt }/{ \cosh^2 t }\,dt = \fr { \pi b}/{ \sinh \pi b},
	\qquad |\Im{b}|<1.
  $$
\PN     %-------------------------------------- IX
\F{\S IX. Rempla\c{c}ant formellement $a$ par $\I\log f$ dans la formule $P$~:}/F
\E{\S IX. Replacing formally $a$ with $\I\log f$ in formula $P$~:}/E
$$
	P(\I\log f,b) =
	\int_0^\infty \fr { \cosh bt }/{ \cosh t + \tfr{(f+\tfr 1/f)}/2 }\,dt = 
		\fr { \pi(f^b-f^{-b}) }/{ (f-f^{-1})\sin\pi b },
		\qquad |\arg f|<\pi,\, |\Re{b}|<1.
 $$
\PN     %-------------------------------------- X
\F{\S X. Autre transform\'ee cosinus de Fourier obtenue en rempla\c{c}ant formellement $b$ par $\I b$ dans la formule pr\'ec\'edente~:}/F
\E{\S X. Another Fourier cosine transform obtained by replacing formally $b$ with $\I b$ in the preceding formula~:}/E
$$
	P(\I\log f,\I b) =
	\int_0^\infty \fr { \cos bt \, dt }/{ \cosh t + \tfr{(f+\tfr 1/f)}/2 } =
			\fr { 2\pi\sin b(\log f) }/{ (f-f^{-1})\sinh \pi b },
			\qquad |\arg f|<\pi,\, |\Im{b}|<1.
 $$
\PN     %-------------------------------------- XIIa
\F{\S XIIa. D\'eveloppement en fractions partielles de la formule $P$ (de la s\'ecante si $\theta=\tfr\pi/2$)~:}/F
\E{\S XIIa. Partial fraction expansion of formula $P$ (of the secant if $\theta=\tfr\pi/2$)~:}/E
$$
	P(\pi-\theta,b) =
	\int_0^\infty \fr { \cosh bt }/{ \cosh t - \cos\theta }\,dt =
		\fr { \pi\sin (\pi-\theta) b }/{ \sin \theta \sin \pi b } =
		\fr 2/{ \sin\theta } \sum_{k=1}^\infty \fr { k\sin k\theta }/{ k^2-b^2 },
		\quad  0<\theta<2\pi,\, |\Re{b}|<1.
 $$
\F{Euler soutient qu'il lui est permis de donner des valeurs purement imaginaires au param\`etre $b$ dans cette derni\`ere s\'erie\footnote{et, implicitement, aussi au param\`etre $b$ de $P$, une m\'ethode employ\'ee dans sa preuve par S.D. Poisson.}, obtenue en partant de la s\'erie r\'ecurrente suivante, qu'il avait d\'emontr\'ee par la m\'ethode des coefficients ind\'etermin\'es\footnote{Leonhard Euler, Inst. Calculi Integralis IV (1845) -- Enestr\"om 660, Supplementum V, \S 187, p. 372.}~:}/F
\E{Euler claims that he is allowed to assign purely imaginary values to $b$ in the last series\footnote{and, implicitly, also to the parameter $b$ of $P$, a method used in its proof by Sim\'eon Denis Poisson.}, obtained by starting from the following recurrent series, which he had proven by the method of undetermined coefficients\footnote{Leonhard Euler, Inst. Calculi Integralis IV -- Enestr\"om 660, Supplementum V, \S 187, p. 372.}~:}/E
$$
	\fr { \sin\theta }/{ 1-2x\cos\theta+x^2 } = \sin\theta + x\sin 2\theta + x^2\sin 3\theta + x^3\sin 4\theta + \etc
 $$
\PN     %-------------------------------------- XIIa
\F{\S XIIa. Un autre d\'eveloppement en fractions partielles (de la s\'ecante hyperbolique si $\theta=\tfr\pi/2$)~:}/F
\E{\S XIIa. Another partial fraction expansion (of the hyperbolic secant if $\theta=\tfr\pi/2$)~:}/E
$$
	P(\pi-\theta,\I b) =
	\int_0^\infty \fr { \cos bt }/{ \cosh t - \cos \theta }\,dt =
	    \fr { \pi\sinh (\pi-\theta) b }/{ \sin \theta\sinh \pi b } =
	    \fr 2/{ \sin \theta } \sum_{k=1}^\infty \fr { k\sin k\theta }/{ k^2+b^2 },
		\quad 0<\theta<2\pi,\, |\Im{b}|<1.
 $$
\PN     %-------------------------------------- XIIb
\F{\S XIIb. Rempla\c{c}ant $t$ par $n t$; $b$ par $\I b$; et $n$ par $\I$ dans la formule $P$ fournit un second paradoxe, soit la valeur purement imaginaire d'une int\'egrale r\'eelle, divergente lorsque prise le long de l'axe r\'eel~:}/F
\E{\S XIIb. Replacing $t$ by $n t$; $b$ by $\I b$; and $n$ by $\I$ in formula $P$ yields a second paradox, a purely imaginary value for a real integral, divergent when taken along the real axis~:}/E
$$
	\int_0^\infty \fr { \cosh bt }/{ \cos t + \cos a }\,dt =
		\fr \pi/{\I} \, \fr { \sinh a b }/{ \sin a\sinh \pi b }\qquad ??
 $$
\F{Poisson a \'etudi\'e ce ph\'enom\`ene\footnote{S.D. Poisson, Sur les int\'egrales des fonctions qui passent par l'infini entre les limites de l'int\'egration, et sur l'usage des Imaginaires dans la d\'etermination des Int\'egrales d\'efinies,  J.\'Ec.Poly. 18/XI(1820) \S 33, 318.}, et a ainsi d\'ecouvert des exemples du th\'eor\`eme des r\'esidus de Cauchy\footnote{Remarque du math\'ematicien fran\c{c}ais Claude Picard, qui a fourni son aide pour la traduction fran\c{c}aise.} en modifiant le chemin d'int\'egration pour \'eviter une singularit\'e.}/F
\E{Poisson investigated this phenomenon\footnote{S.D. Poisson, On the integrals of functions which become infinite between the limits of integration, and on the use of imaginaries in the determination of definite integrals, J. \'Ec. Poly. 18/XI, 1820, \S 33, p. 318.}, and thus discovered examples of the Cauchy residue theorem\footnote{Remark by the French mathematician Claude Picard, who provided help for the French translation.} by changing the path of integration to avoid a singularity.}/E
\medskip
%=================================================================================================
%
\Latex{\end{document}}

\bye